\documentclass[pdflatex]{sn-jnl}
\makeatletter
\let\cl@chapter\undefined
\makeatletter
\usepackage[utf8]{inputenc}
\usepackage{thmtools}
\usepackage{amsmath}
\usepackage{amssymb}
\usepackage{algorithm}
\usepackage{algorithmicx}
\usepackage{graphicx}
\usepackage{color}
\usepackage{a4wide}
\usepackage{hyperref}
\usepackage{cleveref}
\usepackage{tikz}
\usepackage{pgfplots}
\usepackage{todonotes}
\usepackage{algorithm}
\usepackage{algorithmicx}
\usepackage{algpseudocode}

\theoremstyle{thmstyleone}%
\newtheorem{theorem}{Theorem}[section]
\newtheorem{lemma}[theorem]{Lemma}
\newtheorem{corollary}[theorem]{Corollary}
\newtheorem{definition}[theorem]{Definition}
\newtheorem{remark}[theorem]{Remark}

\newtheorem{example}[theorem]{Example}

\newcommand{\dd}{\operatorname{d}\!}
\newcommand{\diam}{\operatorname{diam}}
\newcommand{\dist}{\operatorname{dist}}
\newcommand{\mix}{\operatorname{mix}}

\newcommand{\bfalpha}{\boldsymbol{\alpha}}
\newcommand{\bfbeta}{\boldsymbol{\beta}}
\newcommand{\bfgamma}{\boldsymbol{\gamma}}

\def\letters{a,b,c,d,e,f,g,h,i,j,k,l,m,n,o,p,q,r,s,t,u,v,w,x,y,z,%
	A,B,C,D,E,F,G,H,I,J,K,L,M,N,O,P,Q,R,S,T,U,V,W,X,Y,Z}
\makeatletter
\@for \@l:=\letters \do{%
	\expandafter\edef\csname bb\@l\endcsname{\noexpand\ensuremath{%
  \noexpand\mathbb{\@l}}}%
	\expandafter\edef\csname bf\@l\endcsname{\noexpand\ensuremath{%
	\noexpand\mathbf{\@l}}}%
	\expandafter\edef\csname cal\@l\endcsname{\noexpand\ensuremath{%
  \noexpand\mathcal{\@l}}}%
	\expandafter\edef\csname eu\@l\endcsname{\noexpand\ensuremath{%
  \noexpand\EuScript{\@l}}}%
	\expandafter\edef\csname frak\@l\endcsname{\noexpand\ensuremath{%
  \noexpand\mathfrak{\@l}}}%
	\expandafter\edef\csname rm\@l\endcsname{{\noexpand\rm \@l}}%
	\expandafter\edef\csname scr\@l\endcsname{\noexpand\ensuremath{%
  \noexpand\mathscr{\@l}}}%
}

\DeclareMathOperator{\isdef}{:=}
\DeclareMathOperator{\defis}{=:}

\DeclareMathOperator{\children}{\operatorname{children}}

\DeclareMathOperator{\level}{\operatorname{level}}
\DeclareMathOperator{\depth}{\operatorname{depth}}

\DeclareMathOperator{\OPS}{\operatorname{OPS}}

\newcommand*{\bigtimes}{\mathop{\raisebox{-.5ex}{\hbox{\huge{$\times$}}}}}

\title{Data sparse multilevel covariance estimation in optimal complexity}
\author[1]{J\"urgen D\"olz}\email{doelz@ins.uni-bonn.de}
\affil[1]{\orgdiv{Institute for Numerical Simulation}, \orgname{University of Bonn}, \orgaddress{\street{Friedrich-Hirzebruch-Allee 7}, \city{Bonn}, \postcode{53115}, \country{Germany}}}

\begin{document}
\maketitle
\begin{abstract}%
~We consider the $\calH^2$-formatted compression and computational estimation of covariance functions on a compact set in $\bbR^d$. The classical sample covariance or Monte Carlo estimator is prohibitively expensive for many practically relevant problems, where often approximation spaces with many degrees of freedom and many samples for the estimator are needed. In this article, we propose and analyze a data sparse multilevel sample covariance estimator, i.e., a multilevel Monte Carlo estimator. For this purpose, we generalize the notion of asymptotically smooth kernel functions to a Gevrey type class of kernels for which we derive new variable-order $\calH^2$-approximation rates. These variable-order $\calH^2$-approximations can be considered as a variant of $hp$-approximations. Our multilevel sample covariance estimator then uses an approximate multilevel hierarchy of variable-order $\calH^2$-approximations to compress the sample covariances on each level. The non-nestedness of the different levels makes the reduction to the final estimator nontrivial and we present a suitable algorithm which can handle this task in linear complexity. This allows for a data sparse multilevel estimator of Gevrey covariance kernel functions in the best possible complexity for Monte Carlo type multilevel estimators, which is quadratic. Numerical examples which estimate covariance matrices with tens of billions of entries are presented.

\keywords{covariance estimation \and multilevel Monte Carlo\and $\calH^2$-matrices}
\pacs[MSC Classification]{
65C05 
65D15 
65R20 
}
\end{abstract}

\section{Introduction}
We give a short motivation for fast covariance estimation in \cref{sec:motivation}, a literature review in \cref{sec:literature}, and an outline of the article in \cref{sec:outline}.

\subsection{Motivation}\label{sec:motivation}
Covariance functions or kernel functions
\[
g\colon D\times D\to\bbR,
\]
on a compact set $D\subset\bbR^d$ arise in many fields of application such as Gaussian process computations \cite{RW2006}, machine learning \cite{HSS2008b,Ste2013}, and uncertainty quantification \cite{GHO2017}. However, in many cases these functions are not available in closed form, but must be suitably estimated from samples. Assuming a centered Gaussian random field, the canonical estimator for this purpose is the \emph{sample covariance estimator} or \emph{Monte Carlo estimator}
\[
g\approx\frac{1}{M}\sum_{k=1}^Mz^{(k)}\otimes z^{(k)},
\]
see, e.g., \cite{JW2007}, where the sample functions $z^{(k)}$, $k=1,\ldots,M$, are assumed to be independent, identically distributed (i.i.d.) elements of a Hilbert space with vanishing mean and $\otimes$ is understood as the Hilbertian tensor product. The challenge with the above estimator is that the covariance function and the samples are often infinite-dimensional objects which in practice need to be discretized for computational purposes. After discretization, the sample functions themselves are represented as elements of $\bbR^n$ and the covariance function as a covariance matrix in $\bbR^{n\times n}$.
Assuming that the samples are approximated to an accuracy of $\varepsilon=n^{-\alpha}$, roughly $M=\varepsilon^{-2}=n^{2\alpha}$ samples need to be drawn to reach an overall error of $\calO(\varepsilon)$ of the sample covariance estimator. Thus, the computational effort of the sample covariance estimator is $\calO(Mn^2)=\calO(\varepsilon^{-2-2/\alpha})=\calO(n^{2\alpha+2})$. This is prohibitive for large $n$, as it is often required for sufficient accuracy in applications.

This article presents an algorithm with rigorous error bounds for approximating the covariance function in optimal complexity. Here, optimal complexity is understood such that estimating the covariance has asymptotically the same complexity as estimating the mean, i.e., as good as $\calO(\varepsilon^{-2})=\calO(n^{2\alpha})$ to reach an accuracy of $\calO(\varepsilon)$. To achieve this speedup, we assume a certain behavior, being Gevrey asymptotically smooth, of the covariance function, and that we can draw samples at different levels of resolution from its induced centered Gaussian random field.

\subsection{Related work}\label{sec:literature}
The challenges of large covariance matrices are commonly overcome by using data sparse approximations. Here, the main difference between methods is how the data sparse format is chosen. Purely algebraic methods operate in a black-box fashion on the samples of the sample covariance estimator to estimate suitable compression parameters for previously chosen data sparse formats such as banded matrices \cite{BL2008} or sparse matrices \cite{BL2008a,BL2008,El2008,FLL2016,FGN2006}. See also \cite{CRZ2016} for recent literature review. However, a simultaneous estimate on approximation quality and computational complexity is not available without additional assumptions on the algebraic properties of the samples and/or covariance matrix. These properties are usually inferred from assumed analytical properties of the underlying statistical model.
Here, an often considered analog to some matrix approximation classes considered in \cite{BL2008a} are \emph{asymptotically smooth} covariance functions, which assume a certain decrease of the covariance with increasing spatial distance.
These kinds of functions are also considered in the fast multipole method \cite{GR1987} and it's and abstract counterparts $\calH$- and $\calH^2$-matrices \cite{Bor2010,Hac2015}, as well as in wavelet compression \cite{Sch1998a}. The first have been applied in machine learning \cite{BG2007} and uncertainty quantification \cite{DHS17_1330,HPS2015,KLM2009,ST2006} where complexity and approximation estimates have been derived. The available machinery was also applied to estimate hyperparameters of covariance functions \cite{CGY2016,GAR2019,KLM+2020,LSGK2019,MDHY2017}, but we stress that the objective of this article is to estimate the full covariance functions.
Finally, wavelet based approaches have been used in \cite{HHKS2021,HM2022,HPS2015,HKS2020,SSO2021} for compression and estimation of covariance functions. Similar to wavelet based approaches, sparse grid approaches are also based on a multilevel hierarchy and provide a sparse representation of the covariance matrix, but assume some global smoothness of the covariance \cite{BSZ2011,BG2004a}.
All the mentioned methods operating on assumed analytical properties of covariance functions are capable to reduce the storage requirements of corresponding covariance matrices in $\bbR^{n\times n}$ from $\calO(n^2)$ to $\calO(n)$ or $\calO(n\log^\beta n)$, $\beta>0$, with a negligible approximation error. Thus, the $n^2$ part of the computational cost of the sample covariance estimator can significantly be reduced.

Reducing the computational cost of the sampling process can essentially be achieved by two approaches. The first approach is to see the sample covariance estimator as a Monte Carlo quadrature for a stochastic integral and to replace that quadrature rule by a more efficient method such as quasi-Monte Carlo methods \cite{DKS2013} and sparse grid approaches \cite{BG2004a}. However, bare strong assumptions, further measures to reduce the number of samples are required. The second approach to reduce computational cost during sampling are variance reduction techniques and in particular the multilevel Monte Carlo method, see, e.g., \cite{Gil2015,Hei2001} for a general overview. The basic idea is to exploit a multi-level hierarchy in the approximation spaces for the covariance discretization to obtain covariance matrices of decreasing size and to combine many smaller and only a few larger matrices to a covariance estimator. It was applied to smaller and dense covariance matrices in \cite{MD2019} for the estimation of Sobol indices and to larger covariance matrices combined with a sparse grid approximation in \cite{BSZ2011,CS2013,CS2023} and combined with a wavelet approximation in \cite{HHKS2021}. A multi-fidelity approach was used in \cite{MAPM2025}.

\subsection{Outline}\label{sec:outline}
The article is organized as follows. First, in \cref{sec:results}, we provide an overview over the derived results and introduce the therefore required definitions. \Cref{sec:Gevreyapprox,sec:MC,sec:MLMCalgo} provide omitted proofs. In \ref{sec:experiments} we provide the numerical experiments underlining our theoretical considerations before we draw our conclusions in \cref{sec:concl}.

\section{Statement of results}\label{sec:results}

We give a short overview on the considered covariance function which we aim to estimate in \cref{sec:Gdelta} and discuss the considered approximation schemes in \cref{sec:consideredapprox}. A fast single level sample covariance estimator is presented in \cref{sec:singlelevelintro}, whereas the fast multilevel sample covariance estimator is presented in \cref{sec:multilevelintro}. On the way we introduce notation and preliminary concepts. Proofs are postponed to \cref{sec:Gevreyapprox,sec:MC,sec:MLMCalgo} of this article.

\subsection{Covariance functions considered for estimation}\label{sec:Gdelta}
Let $D\subset\bbR^d$ be compact and equipped with a measure $\mu$. We denote by $L^2(D)=L_\mu^2(D)$ the space of $L^2$-integrable functions with respect to $\mu$ on $D$ and by $H^\gamma(D)$, $\gamma\geq 0$, the classical Sobolev spaces on $D$. We equip $D\times D$ with the product measure $\mu\otimes\mu$. Further, let $(\Omega,\Sigma,\bbP)$ be a probability space and $L_\bbP^p(\Omega;H^\gamma(D))$ the Bochner space of $p$-integrable $H^\gamma(D)$-valued random fields on $\Omega$.

Our aim is to approximate unknown covariance functions $g\colon D\times D\to\bbR$, $g\in H^\gamma(D\times D)$, $\gamma>0$, from samples of the corresponding centered Gaussian random field. We recall that any symmetric, positive semidefinite function $g\colon D\times D\to\bbR$ is the covariance function of a Gaussian random field on $D$ and vice versa. We assume that we know that $g$ is smooth away from the diagonal and satisfies
\begin{align}\label{eq:introgevreysmoothkernel}
|\partial_\bfx^{\bfalpha}\partial_\bfy^{\bfbeta}g(\bfx,\bfy)|\leq C_GA^{|\bfalpha|+|\bfbeta|}(\bfalpha!\bfbeta!)^\delta\|\bfx-\bfy\|_2^{2\theta-d-|\bfalpha|-|\bfbeta|},
\qquad
\bfx,\bfy\in D, \bfx\neq\bfy,
\end{align}
for all $\bfalpha,\bfbeta\in\bbN^d$ and constants $C_G,A>0$,  $\delta\geq 1$, and $\theta$. This is a generalization of the special case $\delta=1$, which is known as asymptotic smoothness \cite{Hac2015}. In analogy to functions of Gevrey class $\delta\geq 1$, being functions  $f\in C^{\infty}(D)$ for which for every $K\Subset D$ and $\bfalpha\in\mathbb{N}^d$ it holds
\[
|\partial^{\bfalpha}f(\bfx)|\leq C_GA^{|\bfalpha|}(\bfalpha!)^\delta
\quad\text{for all}~\bfx\in K,
\]
we refer to \cref{eq:introgevreysmoothkernel} as \emph{Gevrey asymptotic smoothness}, \emph{$G^\delta(C_G,A,2\theta)$-asymptotic smoothness}, or \emph{$G^\delta$-asymptotic smoothness}. See also \cite{CS2012,CVS2011} for a related notion of Gevrey regular singular kernels in the context of ($hp$-)quadrature. Note that we are not requiring any kind of mixed Sobolev smoothness.
\begin{example}
\begin{enumerate}
\item An example for a $G^\delta$-asymptotically smooth covariance function satisfying our assumptions is
\[
g(\bfx,\bfy)
=
\exp\Big(-(2-\|\bfx-\bfy\|_2)^{\frac{1}{1-\delta}}\Big)
\]
with $\bfx,\bfy\in D=\{\bfx\in\bbR^d\colon \|\bfx\|_{2}\leq 1\}$. It can be shown that this function is $G^\delta$-asymptotically smooth, but not asymptotically smooth \cite{CR1996}.

\item More general Gevrey asymptotic smooth random fields can be obtained using stochastic partial differential equation approach to Gaussian random fields \cite{LRL2011,Whi1963}. To this end, let $\calA\colon H^{\theta}(D)\to L^2(D)$ be a self-adjoint and positive definite operator. Denoting by $\calW$ white noise in $L^2(D)$, the solutions to
\[
\calA\calZ=\calW
\]
are realizations of a Gaussian random field $\calZ\in L^2_\bbP(\Omega;H^\theta(D))$. Its covariance operator $\calC=\calA^{-2}$ has explicit representation
\[
(\calC\varphi)(\bfx)=\int_Dg(\bfx,\bfy)\varphi(\bfy)\dd\mu(\bfx),
\]
where $g$ is the covariance function of the Gaussian random field, see \cite[Proposition 2.3]{HHKS2021} for an explicit derivation. For example, the well known Mat\'ern covariance kernels \cite{Mat1960} are given through $D=\bbR^d$ and $\calA=(\kappa^2-\Delta)^{\theta/2}$ with $\kappa>0$, $\theta>d/2$, and are asymptotically smooth. More generally, we may consider any self-adjoint and positive definite pseudo-differential operator $\calA\in\OPS_{cl,\delta}^\theta(D)$ of order $\theta>d/2$ with symbol of Gevrey class $\delta\geq 1$ in the sense of \cite[Definition 1.1]{BK1967}\footnote{We refrain from making this notion more explicit as a proper definition is lengthy and we will not need it for the remainder of the article.}. This implies $\calC=\calA^{-2}\in\OPS_{cl,\delta}^{-2\theta}(D)$ as a consequence of the pseudo-differential operator calculus for Gevrey classes developed in \cite{BK1967}. In analogy to \cite[Lemma 3.0.2]{Sch1998a}
we obtain that the covariance kernel $g$ is $G^\delta$-asymptotically smooth.
\end{enumerate}
\end{example}

\subsection{Considered approximation schemes for single level approximation}\label{sec:consideredapprox}
\subsubsection{Finite element spaces}\label{sec:FEspaces}
For the purpose of our approximation scheme, we assume further that we can generate samples at different levels of resolution from the centered Gaussian random field induced by the unknown covariance function we aim to estimate.
To this end, we consider standard discontinuous, piecewise polynomial finite element spaces $V_h\subset L^2(D)$ of polynomial degree $m$ on a quasi uniform triangulation $\calT_h=\{D_i\}_{i \in I}$. We assume $\calT_h$ to have mesh width $h>0$, and denote the $L^2$-projection onto $V_h$ by $\Pi_h\colon L^2(D)\to V_h$. Denoting by $\otimes$ the Hilbertian tensor product, we identify $L^2(D\times D)\simeq L^2(D)\otimes L^2(D)$ and write ${\Pi_h^{\mix}}=\Pi_h\otimes\Pi_h$ for the $L^2$-projection ${\Pi_h^{\mix}}\colon L^2(D\times D)\to V_h\otimes V_h$. It is well known that the tensorized $L^2$-projection satisfies the approximation estimates.
\begin{align}\label{eq:approximationestimate}
\|g-{\Pi_h^{\mix}}g\|_{L^2(D\times D)}\leq C_{L^2}^\otimes h^{\gamma}\|g\|_{H^{\gamma}(D\times D)},
\qquad
\end{align}
for all $g\in H^\gamma(D\times D)$, $0\leq\gamma\leq m$.

\subsubsection{Cluster trees}\label{sec:ct}
For the purpose of approximation of functions on $D\times D$ we organize the elements of $\calT_h$ into a cluster tree \cite[Chapter 5.3, 5.5, and A.2]{Hac2015} and \cite[Chapter 3.8]{Bor2010}. A \emph{cluster tree} $T_I$ build on the index set $I$ is a tree whose vertices correspond to non-empty subsets of $I$ and are referred to as \emph{clusters}. We require that the root of $T_I$ corresponds to $I$ and that it holds $\dot{\cup}_{s\in\children(t)}s =t$ for all non-leaf clusters $t\in T_I$.
We denote by $L_I$ the leafs of $T_I$ and write 
\[
L_t=\{t_0\in L_I\colon\exists~\text{cluster chain}~t_0\subseteq\ldots\subseteq t_n=t~\text{with}~t_{i-1}\in\children(t_i), i=1,\ldots,n \},
\]
for all $t\in T_I$. The distance of a cluster $t\in T_I$ to the root is denoted by $\level(t)\in\bbN$ and the \emph{depth} of the cluster tree is the maximal level of its clusters. For brevity we write
\begin{equation}\label{eq:depth}
p=\depth(T_I).
\end{equation}
and assume that the cluster trees considered in this article have depth $p=\calO(\log(|I|))$. Moreover, we assume that the number of children for non-leaf clusters is bounded from below and above, i.e.,
\begin{align}\label{eq:ctnoc}
2\leq|\children(t)|\leq C_{\text{ab}},\qquad t\in T_I\setminus L_I,
\end{align}
for some $C_{\text{ab}}>0$, and that the cardinality of the leaf clusters is bounded from below and above, i.e.,
\begin{align}\label{eq:ctleafsize}
n_{\min}/C_{\text{ab}}\leq |t|\leq n_{\min},\qquad t\in L_I,
\end{align}
for some $n_{\min}>0$.
For all $t\in T_I$ we say that $Q_t=\bigtimes_{i=1}^d[a_i,b_i]$ is a \emph{bounding box} of $D_t$ if $D_t=\cup_{i\in t}D_i\subset Q_t$. We assume the bounding boxes $(Q_t)_{t\in T_I}$ to be \emph{$\overline{q}$-regular}. That is, we assume that there is $\overline{q}\in(0,1)$ such that all cluster chains $t_0\subseteq\ldots\subseteq t_n=t$, $t\in T_I$, $t_0\in L_t$, yield families of bounding boxes $(Q_i)_{i=0}^n$, $Q_i=\bigtimes_{j=1}^dJ_j^i$ being bounding boxes to $t_i$, satisfying $|J_j^{i-1}|\leq\overline{q}|J_j^i|$ for all $i=1,\ldots,n$, $j=1,\ldots,d$.

Most standard algorithms for constructing cluster trees from the literature operate on the geometric information of $\calT_h$. When applied to quasi uniform meshes one obtains cluster trees satisfying these assumptions, see also \cite{Bor2010,Hac2015}.

\subsubsection{Block-cluster trees}\label{sec:bct}
All cluster trees can be used to construct a \emph{block-cluster tree} \cite{Bor2010,Hac2015}. A block-cluster tree $T_{I\times I}$ is a tree with vertices corresponding to cluster pairs, referred to as \emph{block-clusters}. Starting with the root $t\times s=I\times I$, the block-clusters of the tree are checked recursively whether they fulfill the admissibility condition
\begin{align}\label{eq:admissibility}
\max\big\{\diam_\infty Q_t,\diam_\infty Q_s\}=\diam_\infty(Q_t\times Q_s)\leq 2\eta\dist_2(Q_t,Q_s).
\end{align}
When admissible, they are added to the set of \emph{far-field} block clusters $L_{I\times I}^+$. Otherwise, the recursion is applied to all block clusters $t'\times s'$, $t'\in\children(t)$, $s'\in\children(s)$. If $t$ or $s$ have no children $t\times s$ is added to $L_{I\times I}^-$. Clearly, this generates a tree structure on $T_{I\times I}$ whose set of leafs is given as $L_{I\times I}=L_{I\times I}^+\cup L_{I\times I}^-$. For later reference we write
\[
L_{t\times s}=\{t_0\times s_0\colon t_0\in L_t, s_0\in L_s \}.
\]

\subsubsection{Variable-order $\calH^2$-space}
With preliminaries in place, we are in the position to introduce the approximation space for Gevrey asymptotically smooth functions. To this end, let $\alpha\in\bbN_0$, $\beta\in\bbN$,
\begin{equation}\label{eq:kidelta}
k_i^\delta=\lceil(\beta+\alpha i)^\delta\rceil,\qquad i\in\bbN.
\end{equation}
Denote by $\calI_k^{Q}$ the tensorized Chebyshev interpolation operator on the bounding box $Q$ with $k+1$ interpolation points per direction. For $t,s\in T_I$, $t_0\in L_t$, $s_0\in L_s$ and $t_0\subseteq\ldots\subseteq t_n=t$ and $s_0\subseteq\ldots\subseteq s_m=s$ cluster chains in $T_I$, we define the interpolation operators
\begin{equation}\label{eq:Ir}
\calI_{t_0}^{t}=\calI_{t_0}\circ\ldots\circ\calI_{t_n},
\qquad
\text{with}~\calI_{t_i}=\calI_{k_{p-\level(t_i)}^\delta}^{Q_i}~\text{for}~i=0,\ldots,n,
\end{equation}
with $p$ as in \cref{eq:depth}, and
\[
\calI_{t_0\times s_0}^{t\times s}=\calI_{t_0}^{t}\otimes\calI_{s_0}^{s}.
\]
We define $\calP_k$ as the space of polynomials on $\bbR^d$ with degrees $\|\bfalpha\|_\infty\leq k$ and set
\begin{equation}\label{eq:Pts}
\calP_{t\times s} = \big(\calP_{k_{p-\level(t)}^\delta}\otimes\calP_{k_{p-\level(s)}^\delta}\big)\big|_{D_t\times D_s}
\end{equation}
for all $t,s\in T_I$ and
\begin{equation}\label{eq:Ptspw}
\calP_{t\times s}^{\text{\text{pw}}} = \{f\colon D_t\times D_s\to\bbR\colon f=\calI_{t_0\times s_0}^{t\times s}q,t_0\times s_0\in L_{t\times s},q\in\calP_{t\times s}\} 
\end{equation}
for all $t\times s\in L_{I\times I}^+$.

Finally, we define the \emph{$\calH^2$-space of kernel functions} as
\begin{equation}\label{eq:H2space}
V^\calH=\Big\{g\colon D\times D\to\bbR\colon g\big|_{D_t\times D_s}\in\calP_{t\times s}^{\text{pw}}~\text{for all}~t\times s\in L_{I\times I}^+\Big\}.
\end{equation}
For later reference we note that the spaces $V^\calH$ use
\begin{align}\label{eq:rankdistribution}
K_t=\big(k_{p-\level(t)}^\delta\big)^d=\big\lceil(\beta+\alpha(p-\level(t)))^\delta\big\rceil^{d}
\end{align}
degrees of freedom per cluster $t\in T_I$.

As we will discuss in \cref{lem:H2complexity} below, the significance of this function space is that only $\calO(|I|)$ function values of $g$ are required to represent $g|_{L_{I\times I}^+}$. We also note that $V^\calH$ has a vector space structure. For implementation purposes, we note that $V^\calH$ does not need to be implemented, but naturally arises from parameter choices in existing code bases.

\subsubsection{Approximation of Gevrey kernels in variable-order $\calH^2$-spaces}\label{sec:approxGvoH2}
Let
\begin{equation}\label{eq:defPih}
\Pi^\calH\colon L^2(D\times D)\to V^\calH
\end{equation}
denote the $L^2(D\times D)$-projection onto $V^\calH$. Assuming that there are constants $C_{\text{cu}}\in\bbR_{>0}$, $h_{\calH}\in\bbR_{>0}$, $C_{\text{gr}}\in\bbR_{>0}$, and $\zeta\in\bbR_{\geq 1}$ such that
\[
\mu(D_t)\leq C_{\text{cu}}\diam_{\infty}(Q_t)^d, 
\]
for all $t\in T_I$,
\[
C_{\text{gr}}^{-1}h_{\calH}\leq\diam_{\infty}(Q_t)\leq C_{\text{gr}}h_{\calH}
\]
for all $t\in L_I$, and
\[
\diam_{\infty}(Q_t)\leq\zeta\diam_{\infty}(Q_{t'})
\]
for all $t'\in\children(t)$, $t\in T_I$, see also \cite[(4.58) and (4.59)]{Bor2010}, we show in \cref{thm:L2interpolationestimate} approximation estimates for $\Pi^\calH$. For ease of presentation we state in the following a simplified version.
\begin{theorem}[{Simplified version of \cref{thm:L2interpolationestimate}}]\label{thm:kernell2errorintro}
There are $\tilde{\rho}\in(0,1)$ (specified in \cref{eq:rhotilde}), $\alpha\in\bbN$, and $\beta_0\in\bbN$ such that $\zeta^{-2\theta}\tilde{\rho}^{\alpha}<1$ and
\[
\big\|g-\Pi^\calH g\big\|_{L^2(D\times D)}
\leq
C\frac{h_{\calH}^{-2\theta}\tilde{\rho}^{\beta}}{1-\zeta^{-2\theta}\tilde{\rho}^{\alpha}}
\]
for a fixed constant $C>0$ and for all $\beta\geq\beta_0$.
\end{theorem}
Thus, as immediate from the theorem, $G^\delta$-asymptotically smooth functions can be arbitrarily well approximated in the $\calH^2$-format, with exponential convergence rates in $\beta$. We recall that $\beta$ is the minimal polynomial degree for polynomial interpolation, cf.\ \cref{eq:kidelta}.

This approximation result is based on a new approximation result for Gevrey functions in one dimension, which we state in the following for the reader's convenience and prove in \cref{cor:abGevreyInterpolation} below. In comparison to other approximation results in the literature, we note that the dependence of the contraction factor on $A$ is explicit and convergence is guaranteed independently of its magnitude. This is an essential ingredient for establishing the $\calH^2$-approximation rates, independently of the chosen admissibility parameter $\eta$ in \cref{eq:admissibility}.

\begin{theorem}[{\Cref{cor:abGevreyInterpolation}}]\label{thm:abGevreyInterpolationintro}
Let $\rho(r)=r+\sqrt{1+r^2}$. For any $f\in G^\delta([a,b],C_G,A)$, $B=A(b-a)/2$, and $m\in\bbN$, $m\geq 3$, it holds that
\[
\min_{q\in\calP_m}\|f-q\|_{C([a,b])}\leq C(B,\delta)C_G\rho(1/B)^{-m^{1/\delta}/e^2},
\]
where $C(B,\delta)$ is monotonically increasing in $B$.
\end{theorem}

For $\delta=1$, i.e., for analytic functions, the approximation results imply exponential convergence, similarly to other results in the literature. On the other hand, with increasing $\delta$, the approximation rates reduce due to decreasing regularity.

\subsection{$\calH^2$-formatted sample covariance estimator}\label{sec:singlelevelintro}

\subsubsection{The estimator}
For the purpose of a fast single level sample covariance estimator we assume that we can draw iid discretized samples $\Pi_hz^{(k)}\in V_h$, $k=1,2,\ldots$ from the centered Gaussian random field induced by the covariance function. More precisely, using the notation of the previous subsections, the (single level) $\calH^2$-formatted sample covariance estimator is given as follows.
\begin{definition}\label{def:H2MCFEintro}
Given discretized iid samples $\Pi_hz^{(k)}$, $k=1,\ldots,M$, of the centered Gaussian random field induced by the covariance function $g$, the \emph{(single level) $\calH^2$-formatted sample covariance estimator ($\calH^2$-SCE)} is defined as
\begin{equation}\label{eq:H2SCE}\tag{$\calH^2$-SCE}
\bbE[\Pi^\calH\Pi_h^{\mix}g]
\approx E^{MC}[\Pi^\calH\Pi_h^{\mix}g]
=
\frac{1}{M}\sum_{k=1}^M\Pi^\calH\Big(\Pi_hz^{(k)}\otimes \Pi_hz^{(k)}\Big).
\end{equation}
\end{definition}

Given the discretized samples $\Pi_hz^{(k)}\in V_h$, the $\calH^2$-SCE relies on two main algorithmic operations: The $L^2$-projection $\Pi^\calH$ of the dyads $\Pi_hz^{(k)}\otimes \Pi_hz^{(k)}$ into $V^\calH$ and the summation over $k$. Due to the vector space structure of $V^\calH$, computing the sum is straightforward, and, after discretization, corresponds to the $\calH^2$-matrix summation. This operation is known to have linear complexity \cite[Chapter 7.3]{Bor2010}. Thus, given $\Pi_hz^{(k)}$, it remains to clarify how $\Pi^\calH(\Pi_hz^{(k)}\otimes \Pi_hz^{(k)})$ can be efficiently implemented. 

In the following, we present an algorithm which allows computing $\Pi^\calH\Big(\Pi_hz^{(k)}\otimes \Pi_hz^{(k)}\Big)$ in $\calO(|I|)$ complexity. Thus, the overall complexity of the $\calH^2$-SCE is $\calO(M|I|)$.

\subsubsection{Algorithmic realization}
For the algorithmic implementation of $\Pi^\calH(\Pi_hz^{(k)}\otimes \Pi_hz^{(k)})$ we consider for all $r\in L_I$ the \emph{moment matrices} $\bfM_r=[(\psi_i^r,\phi_j^r)_{L^2(D_r)}]_{i,j}$, where $\calP_r=\operatorname{span}\{\psi_i^r\}_{i=1}^{m}$ is a nodal basis of $\calP_r$, and $V_h|_{D_r}=\operatorname{span}\{\phi_i^r\}_{i}$ is the used finite element basis of $V_h$ restricted to $r$. Moreover, for all $r\in T_I\setminus L_I$, $r'\in\children(r)$, we consider the \emph{transfer matrices} $\bfE_{r'}$. These matrices are the matrix representation of $\calE_{r'}\colon\calP_r\to\calP_{r'}$ defined by the interpolation $q\mapsto\calI_{r'}q$, cf.\ \cref{eq:Ir}, with respect to the nodal bases $\{\psi_i^r\}_{i=1}^m$ and $\{\psi_i^{r'}\}_{i=1}^m$. For the constant order case, i.e., for $\alpha=0$, we denote the family of transfer matrices by $\{\bfF_t\}_{t\in T_I\setminus\{I\}}$. These matrices are standard in the context of $\calH^2$-matrices \cite{Bor2010}. For the purpose of implementing the $\calH^2$-SCE, we further require the matrices $\bfQ_r=\big[(\psi_i^r,\psi_j^r)_{L^2(D_r)}\big]_{i,j=1}^m$. Given these preliminaries, an efficient algorithm to compute $\Pi^{\calH}(\Pi_hz\otimes\Pi_hz)\in\Pi^{\calH}(V_h\otimes V_h)$ is given in \cref{alg:Pih}. A derivation of this algorithm is given in \cref{sec:Pih}.
\begin{algorithm}
\caption{\label{alg:Pih}Compute $\Pi^\calH(\Pi_hz\otimes\Pi_hz)$ from coefficient vector $\bfz^h$ of $\Pi_hz$}
\begin{algorithmic}[1]
\State \textbf{Input:} Moment matrices $\{\bfM_r\}_{r\in L_I}$, transfer matrices $\{\bfE_r\}_{r\in T_I}$, matrices $\{\bfQ_r\}_{r\in L_I}$
\State \phantom{\textbf{Input:} }coefficient vector $\bfz^h$ of $\Pi_hz\in V_h$
\State \textbf{Output:} Coefficients $\{\bfU_{t\times s}\}_{t\times s\in L_{I\times I}}$ of $\Pi^{\calH}(\Pi_hz\otimes\Pi_hz)\in\Pi^{\calH}(V_h\otimes V_h)$.
\bigskip
\State Compute $\bfq_r^h=\bfM_r\bfz_r^h$ for all $r\in L_I$.
\State Compute $\bfq_r^h=\sum_{r'\in\children(r)}\bfE_{r'}^\intercal \bfq_{r'}^h$ for all $r\in T_I\setminus L_I$ recursively.
\State Compute $\bfQ_r=\sum_{r'\in\children(r)}\bfE_{r'}^\intercal \bfQ_{r'}\bfE_{r'}$ for all $r\in T_I\setminus L_I$ recursively.
\State Solve the local systems $\bfQ_r\bfu_r=\bfq_r^h$ for all $r\in T_I$.
\State Set $\bfU_{t\times s}=\bfu_t\otimes\bfu_s$ for all $t\times s\in L_{I\times I}^+$ and $\bfU_{t\times s}=\bfz_t^h\otimes\bfz_s^h$ for all $t\times s\in L_{I\times I}^-$
\end{algorithmic}
\end{algorithm}

The following remarks indicate that the implementation efforts for \cref{alg:Pih} are rather minor.
\begin{remark}\label{rem:algorem}
\begin{itemize}
\item The moment matrices $\{\bfM_r\}_{r\in L_I}$ and the transfer matrices $\{\bfE_r\}_{r\in T_I}$ are implemented in most $\calH^2$-codes and can be computed in linear complexity in $|I|$, see \cite{Bor2010}. The matrices $\{\bfQ_r\}_{r\in L_I}$ amount to $L^2$-inner products of polynomials, and are thus easily implemented. Following similar arguments, they can be computed in linear complexity in $|I|$ as well.
\item Line 5 in the algorithm is known as \emph{forward transformation}, and also implemented in most $\calH^2$-codes.
\item Exploiting that $\operatorname{vec}(\bfE_{r'}^\intercal \bfQ_{r'}\bfE_{r'})=(\bfE_{r'}\otimes \bfE_{r'})\operatorname{vec}(\bfQ_{r'})$, with $\operatorname{vec}(\cdot)$ being the reshaping of the columns of matrix to a long column vector, line 6 can also be implemented using the forward transformation.
\item In actual implementations we may compute and factorize $\{\bfQ_t\}_{t\in T_I}$ once and use it for all samples, whereas $\{\bfq_t\}_{t\in T_I}$ needs to be recomputed for each sample.
\end{itemize}
\end{remark}

We show in \cref{sec:Pihcompl} that \cref{alg:Pih} has complexity $\calO((\alpha+\beta)^{2\delta d}|I|)$, if the samples $\Pi_hz^{(k)}$ are available. An algorithm for the $\calH^2$-SCE is given in \cref{alg:H2SCE}.

\begin{algorithm}
\caption{\label{alg:H2SCE}Algorithmic implementation of the $\calH^2$-SCE}
\begin{algorithmic}[1]
\State \textbf{Input:} Finite element space $V_h$, number of samples $M$, $\alpha$, $\beta$, $\delta$
\State \textbf{Output:} Coefficients $\{\bfU_{t\times s}\}_{t\times s\in L_{I\times I}}$ of the $\calH^2$-SCE
\bigskip
\State Compute cluster tree $T_I$ and block cluster tree $T_{I\times I}$ using standard codes
\State Compute moment matrices $\{\bfM_r\}_{r\in L_I}$, transfer matrices $\{\bfE_r\}_{r\in T_I}$ using standard codes
\State Compute matrices $\{\bfQ_r\}_{r\in L_I}$
\smallskip
\State Initialize $\{\bfU_{t\times s}\}_{t\times s\in L_{I\times I}}$ to zero
\For{$k=1,\ldots,M$}
\State Draw sample vector $\bfz_k^h$ representing $\Pi_hz^{(k)}$
\State Compute coefficients $\{\bfU_{t\times s}^{(k)}\}_{t\times s\in L_{I\times I}}$ of $\Pi^\calH(\Pi_hz^{(k)}\otimes\Pi_hz^{(k)})$ using \cref{alg:Pih}
\State Update $\{\bfU_{t\times s}\}_{t\times s\in L_{I\times I}}+\!\!=\{\bfU_{t\times s}^{(k)}\}_{t\times s\in L_{I\times I}}$ using $\calH^2$-matrix addition from standard codes
\EndFor
\end{algorithmic}
\end{algorithm}

As indicated above, the complexity of \cref{alg:H2SCE} is $\calO(M(\alpha+\beta)^{2\delta d}|I|)$.

\subsubsection{Approximation result}

We obtain the following approximation result for the $\calH^2$-SCE, which we show later in \cref{lem:covHL2error}.

\begin{theorem}[{Simplified version of \cref{lem:covHL2error}}]\label{lem:covHL2errorintro}
Let $2\theta<0$ and let $g\in H^\gamma(D\times D)$ be $G^\delta(C_G,A,2\theta)$-asymptotically smooth. Then there are $\tilde{\rho}\in(0,1)$ (specified in \cref{eq:rhotilde}), $\alpha\in\bbN$, and $\beta_0\in\bbN$ such that $\zeta^{-2\theta}\tilde{\rho}^{\alpha}<1$ and
\[
\big\|g-E^{MC}[\Pi^\calH\Pi_h^{\mix}g]\big\|_{L_{\bbP}^2(\Omega,L^2(D\times D))}
\leq
\frac{C_{\text{lc}}C_{\text{sp}}h_{\calH}^{-2\theta}\tilde{\rho}^{\beta}}{1-\zeta^{-2\theta}\tilde{\rho}^{\alpha}}+\bigg(C_{L^2}^\otimes h^{\gamma}+ \frac{1}{\sqrt{M}}\bigg)\|g\|_{H^\gamma(D\times D)}
\]
for all $\beta\geq\beta_0$.
\end{theorem}

Thus, by increasing $\beta$, $M$, and by decreasing $h$, we can approximate the covariance function $g$ arbitrarily well with the $\calH^2$-SCE. The only disadvantage is the overall complexity of
\begin{equation}\label{eq:H2SCEcompl}
\calO(M(\alpha+\beta)^{2\delta d}|I|)
=
\calO(M(\alpha+\beta)^{2\delta d}h^{-d}).
\end{equation}
This is prohibitive for large $M$ and small $h$ and will be addressed in \cref{sec:multilevelintro}.

\subsubsection{Computational $\calH^2$-sample covariance estimation}
For computational covariance estimation one often aims at a discretization of the covariance function rather than the covariance itself. In the following we provide error estimates for bilinear forms of type
\begin{align}\label{eq:galerkinblf}
	a(u_h,v_h)=\int_D\int_Dg(\bfx,\bfy)u_h(\bfx)v_h(\bfy)\dd\mu(\bfx)\dd\mu(\bfy)
\end{align}
for $u_h,v_h\in V_h$. Of course, any other discontinuous finite element space other than $V_h$ on $\calT$ could also be used. Using the Cauchy-Schwartz estimate \cref{lem:covHL2errorintro} implies
\begin{corollary}\label{lem:fembilinearerror}
	Let the assumptions of \cref{lem:covHL2error} hold. Then it holds
	\begin{align*}
		&\bigg\|\int_{D}\int_{D}\Big(g(\bfx,\bfy)-E^{MC}[\Pi^{\calH}\Pi_h^{\mix}g(\bfx,\bfy)]\Big)u_h(\bfx)v_h(\bfy)\dd\mu(\bfx)\dd\mu(\bfy)\bigg\|_{L_{\bbP}^2(\Omega)}\\
		&\qquad\leq \bigg(\frac{C_{\text{lc}}C_{\text{sp}}h_{\calH}^{-2\theta}\tilde{\rho}^\beta}{1-\zeta^{-2\theta}\tilde{\rho}^\alpha}+\bigg(C_{L^2}^\otimes h^{{\gamma}}+ \frac{1}{\sqrt{M}}\bigg)\|g\|_{H^\gamma(D\times D)}\bigg)\|u_h\|_{L^2(D)}\|v_h\|_{L^2(D)},
	\end{align*}
	for all $u_h,v_h\in V_h$ and $\beta\geq\beta_0$.
\end{corollary}
For a finite element basis $V_h=\operatorname{span}\{\phi_i\}_{i\in I}$ our assumptions guarantee that $\bfA=[a(\phi_j,\phi_i)]_{i,j\in I}$ can be stored as $\calH^2$-matrix with a storage requirement of $C_{\calH^2}(\alpha+\beta)^{\delta d}|I|$, i.e., linear in the cardinality of $I$, cf.\ \cite{Bor2010}. Thus, using \cref{alg:Pih}, the $\calH^2$-SCE is computable in \cref{eq:H2SCEcompl} complexity also when bilinear forms of the type \cref{eq:galerkinblf} are considered.

We remark that methods relying on a sparse grid approximation of the covariance  yield a complexity which is only linear up to a logarithmic factor, see, e.g., \cite{BSZ2011}.

To overcome the computational prohibitive cost \cref{eq:galerkinblf} of the $\calH^2$-SCE for large $M$ and small $h$, we derive in the following a multilevel version thereof.

\subsection{Multilevel $\calH^2$-formatted sample covariance estimator}\label{sec:multilevelintro}

\subsubsection{Construction of multilevel hierarchies}\label{sec:mlhier}
\begin{figure}
	\begin{tikzpicture}[scale=7.5]
		\draw (0,0.065) node {\begin{tikzpicture}[level distance=10mm,
every node/.style={fill=red!60,rectangle,rounded corners,inner sep=1pt},
level 1/.style={sibling distance=40mm,nodes={fill=red!50}},
level 2/.style={sibling distance=20mm,nodes={fill=red!40}},
level 3/.style={sibling distance=6.5mm,nodes={fill=red!30}}]
\node {$I_0=\{1,2,3\}$}
	child {
		child {node {$\{1\}$}}
	}
	child {node {$\{2,3\}$}
		child {node {$\{2\}$}}
		child {node {$\{3\}$}}
	};
\end{tikzpicture}};
		\draw (1,0) node {\begin{tikzpicture}[level distance=10mm,
every node/.style={fill=red!60,rectangle,rounded corners,inner sep=1pt},
level 1/.style={sibling distance=40mm,nodes={fill=red!50}},
level 2/.style={sibling distance=20mm,nodes={fill=red!40}},
level 3/.style={sibling distance=6.5mm,nodes={fill=red!30}},
level 3/.style={sibling distance=6.5mm,nodes={fill=red!20}}]
\node {$I_1=\{1,\ldots,9\}$}
	child {child {node {$\{1,2,3\}$}
		child {node {$\{1\}$}}
		child {node {$\{2\}$}}
		child {node {$\{3\}$}}
	}}
	child {node {$\{4,\ldots,9\}$}
		child {node {$\{4,5,6\}$}
			child {node {$\{4\}$}}
			child {node {$\{5\}$}}
			child {node {$\{6\}$}}
		}
		child {node {$\{7,8,9\}$}
			child {node {$\{7\}$}}
			child {node {$\{8\}$}}
			child {node {$\{9\}$}}
		}
	};
\end{tikzpicture}};
		\draw (0,-0.4) node {\begin{tikzpicture}[scale=2.2]
\draw (0,0) -- (3,0);
\foreach \i in {0,...,3}
{
\draw (\i,-0.1) -- (\i,0.1);
}
\foreach \i in {1,...,3}
{
	\draw (\i-0.5+0.02,0.15) node {$D_{\i}^{(0)}$};
}
\end{tikzpicture}};
		\draw (1,-0.4) node {\begin{tikzpicture}[scale=2.2]
\draw (0,0) -- (3,0);
\foreach \i in {0,...,3}
{
\draw (\i,-0.1) -- (\i,0.1);
}
\foreach \i in {0,...,9}
{
	\draw (\i/3,-0.05) -- (\i/3,0.05);
}
\foreach \i in {1,...,9}
{
	\draw (\i/3-1/6+0.02,0.15) node {$D_{\i}^{(1)}$};
}
\end{tikzpicture}};
	\end{tikzpicture}
	\caption{\label{fig:nestedclustertree}Illustration of ``nested'' cluster trees $T_{I_0}$ (upper left) and $T_{I_1}$ (upper right) to nested decompositions $\{D_i^{(0)}\}_{i\in I_0}$ (bottom left) and $\{D_i^{(1)}\}_{i\in I_1}$ (bottom right).}
\end{figure}
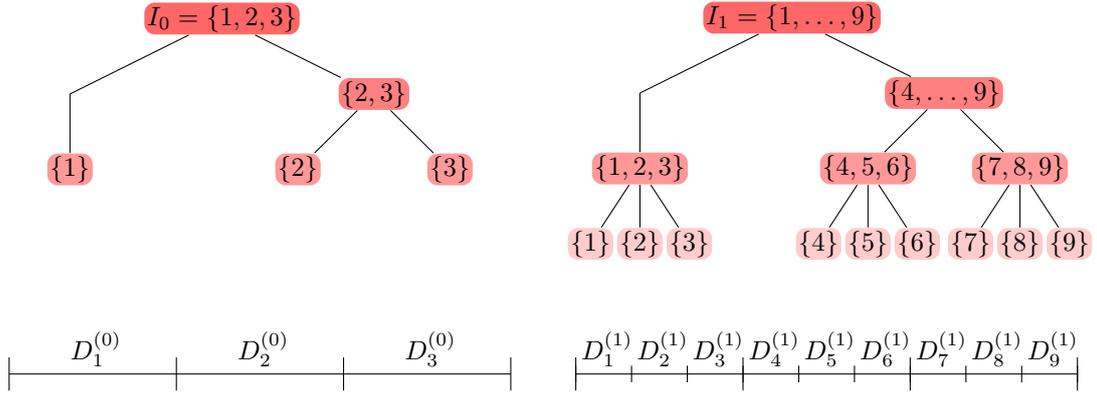

We start with a discontinuous, piecewise polynomial finite element space $V_{h_0}\subset L^2(D)$ of degree $m$ on a quasi uniform triangulation $\mathcal{T}_{h_0}=\{D_i^{(0)}\}_{i\in I_0}$, as in \cref{sec:FEspaces}. We further assume to have a cluster tree $T_{I_0}$ be a cluster tree constructed on $I_0$ to our disposal, as in \cref{sec:ct}. Under these circumstances we can generate a sequence of nested decompositions
\begin{align}\label{eq:Cuni}
	\Big\{\mathcal{T}_{h_\ell}=\big\{D_i^{(\ell)}\big\}_{i\in I_\ell}\Big\}_{\ell=0}^\infty
	\qquad\text{with}\qquad
	|I_\ell|=|I_0|C_{\text{uni}}^{\ell}
\end{align}
for some $C_{\text{uni}}> 1$ and corresponding finite element spaces $V_{h_0}\subset V_{h_1}\subset V_{h_2}\subset\ldots\subset L^2(D)$ in the usual way using uniform refinement. In \cref{alg:nestedct} we provide an algorithm which, called recursively, constructs ``nested'' cluster trees on the triangulations.
\begin{algorithm}
\caption{\label{alg:nestedct}Construction of ``nested'' cluster trees}
\begin{algorithmic}[1]
\State \textbf{Input:} $\calT_{h_\ell}=\{D_i^{(\ell)}\}_{i\in I_\ell}$, $T_{I_\ell}$ cluster tree on $\calT_{h_\ell}$, $\calT_{h_{\ell+1}}=\{D_i^{(\ell+1)}\}_{i\in I_{\ell+1}}$
\State \textbf{Output:} Cluster tree $T_{I_{\ell+1}}$ on $\calT_{h_{\ell+1}}$
\bigskip
\State Set $T_{I_{\ell+1}}=T_{I_\ell}$
\State Replace all $t\in T_{I_{\ell+1}}$ by their uniform refinement
\State For all $t\in L_{I_{\ell+1}}$, redefine $\children(t)=\{D_i^{(\ell)}\}_{i\in t}$
\end{algorithmic}
\end{algorithm}
For an illustration of those ``nested'' cluster trees see \cref{fig:nestedclustertree}.
In the following, we will assume that we have cluster trees $T_{I_{\ell+1}}$, $\ell=1,2,\ldots$, constructed from \cref{alg:nestedct} to our disposal. These generated cluster trees satisfy the assumptions from \cref{sec:ct}.

The nestedness of the generated cluster trees directly implies that also the sequence of block-cluster trees $\{T_{I_\ell\times I_\ell}\}_{\ell=1}^{\infty}$ constructed as in \cref{sec:bct} is ``nested''. To simplify notation we introduce the equivalence relation
\[
	t_\ell\sim t_{\ell+1}\qquad\text{for all}\qquad t_\ell\in T_{I_\ell},~t_{\ell+1}\in T_{I_{\ell+1}}~\text{for which}~\cup_{k\in t_\ell}D_k^{(\ell)}=\cup_{k\in t_{\ell+1}}D_k^{(\ell+1)},
\]
and similarly for the block-cluster tree. We refer to clusters up to equivalence relations, meaning that clusters can be elements of several cluster trees. Likewise, block-clusters can be elements of several block-cluster trees.

In the following we exploit that the block-cluster trees generated from ``nested'' cluster trees are ``nested'' as well. Thus, they provide a sequence of nested decompositions of $D\times D$. Note that the farfields and the nearfields of ``nested'' block-cluster trees do not provide nested decompositions of $D\times D$, since only
\begin{equation}\label{eq:adminimpladmin}
t\times s\in L^+_{I_\ell\times I_\ell}\Rightarrow t\times s\in L^+_{I_{\ell+1}\times I_{\ell+1}}
\end{equation}
is guaranteed from the construction. Thus, the sequence $\{V^{\calH_\ell}\}_{\ell=0}^\infty$ of $\calH^2$-spaces from \cref{eq:H2space} generated by the sequence of block-cluster trees is not nested, see also \cref{fig:nestedVh} for an illustration.
\begin{figure}
	\centering
	\begin{tikzpicture}[scale=0.13]

\foreach \x in {0,...,3} {
	\draw[fill=red!50] (8*\x,24-8*\x) rectangle (8+8*\x,32-8*\x);
}
\foreach \x in {0,...,2} {
	\draw[fill=red!50] (8*\x,16-8*\x) rectangle (8+8*\x,24-8*\x);
	\draw[fill=red!50] (8+8*\x,24-8*\x) rectangle (16+8*\x,32-8*\x);
}

\foreach \x in {0,...,1} {
	\draw[fill=blue!20] (8*\x,8-8*\x) rectangle (8+8*\x,16-8*\x);
	\draw[fill=blue!20] (16+8*\x,24-8*\x) rectangle (24+8*\x,32-8*\x);
	\draw (4+8*\x,12-8*\x) node {\large 9};
	\draw (20+8*\x,28-8*\x) node {\large 9};
}
\draw[fill=blue!20] (0,0) rectangle (8,8);
\draw[fill=blue!20] (24,24) rectangle (32,32);
\draw (4,4) node {\large 9};
\draw (28,28) node {\large 9};

\end{tikzpicture}
	\qquad
	\begin{tikzpicture}[scale=0.13]

\foreach \x in {0,...,7} {
	\draw[fill=red!50] (4*\x,28-4*\x) rectangle (4+4*\x,32-4*\x);
}
\foreach \x in {0,...,6} {
	\draw[fill=red!50] (4*\x,24-4*\x) rectangle (4+4*\x,28-4*\x);
	\draw[fill=red!50] (4+4*\x,28-4*\x) rectangle (8+4*\x,32-4*\x);
}

\foreach \x in {0,...,5} {
	\draw[fill=blue!10] (4*\x,20-4*\x) rectangle (4+4*\x,24-4*\x);
	\draw[fill=blue!10] (8+4*\x,28-4*\x) rectangle (12+4*\x,32-4*\x);
	\draw (2+4*\x,22-4*\x) node {\small 9};
	\draw (10+4*\x,30-4*\x) node {\small 9};
}
\foreach \x in {0,...,2} {
	\draw[fill=blue!10] (8*\x,16-8*\x) rectangle (4+8*\x,20-8*\x);
	\draw[fill=blue!10] (12+8*\x,28-8*\x) rectangle (16+8*\x,32-8*\x);
	\draw (2+8*\x,18-8*\x) node {\small 9};
	\draw (14+8*\x,30-8*\x) node {\small 9};
}

\foreach \x in {0,...,1} {
	\draw[fill=blue!20] (8*\x,8-8*\x) rectangle (8+8*\x,16-8*\x);
	\draw[fill=blue!20] (16+8*\x,24-8*\x) rectangle (24+8*\x,32-8*\x);
	\draw (4+8*\x,12-8*\x) node {\large 25};
	\draw (20+8*\x,28-8*\x) node {\large 25};
}
\draw[fill=blue!20] (0,0) rectangle (8,8);
\draw[fill=blue!20] (24,24) rectangle (32,32);
\draw (4,4) node {\large 25};
\draw (28,28) node {\large 25};

\end{tikzpicture}
	\qquad
	\begin{tikzpicture}[scale=0.13]

\foreach \x in {0,...,15} {
	\draw[fill=red!50] (2*\x,30-2*\x) rectangle (2+2*\x,32-2*\x);
}
\foreach \x in {0,...,14} {
	\draw[fill=red!50] (2*\x,28-2*\x) rectangle (2+2*\x,30-2*\x);
	\draw[fill=red!50] (2+2*\x,30-2*\x) rectangle (4+2*\x,32-2*\x);
}

\foreach \x in {0,...,13} {
	\draw[fill=blue!10] (2*\x,26-2*\x) rectangle (2+2*\x,28-2*\x);
	\draw[fill=blue!10] (4+2*\x,30-2*\x) rectangle (6+2*\x,32-2*\x);
	\draw (1+2*\x,27-2*\x) node {\tiny 9};
	\draw (5+2*\x,31-2*\x) node {\tiny 9};
}
\foreach \x in {0,...,6} {
	\draw[fill=blue!10] (4*\x,24-4*\x) rectangle (2+4*\x,26-4*\x);
	\draw[fill=blue!10] (6+4*\x,30-4*\x) rectangle (8+4*\x,32-4*\x);
	\draw (1+4*\x,25-4*\x) node {\tiny 9};
	\draw (7+4*\x,31-4*\x) node {\tiny 9};
}

\foreach \x in {0,...,5} {
	\draw[fill=blue!20] (4*\x,20-4*\x) rectangle (4+4*\x,24-4*\x);
	\draw[fill=blue!20] (8+4*\x,28-4*\x) rectangle (12+4*\x,32-4*\x);
	\draw (2+4*\x,22-4*\x) node {\small 25};
	\draw (10+4*\x,30-4*\x) node {\small 25};
}
\foreach \x in {0,...,2} {
	\draw[fill=blue!20] (8*\x,16-8*\x) rectangle (4+8*\x,20-8*\x);
	\draw[fill=blue!20] (12+8*\x,28-8*\x) rectangle (16+8*\x,32-8*\x);
	\draw (2+8*\x,18-8*\x) node {\small 25};
	\draw (14+8*\x,30-8*\x) node {\small 25};
}

\foreach \x in {0,...,1} {
	\draw[fill=blue!30] (8*\x,8-8*\x) rectangle (8+8*\x,16-8*\x);
	\draw[fill=blue!30] (16+8*\x,24-8*\x) rectangle (24+8*\x,32-8*\x);
	\draw (4+8*\x,12-8*\x) node {\large 49};
	\draw (20+8*\x,28-8*\x) node {\large 49};
}
\draw[fill=blue!30] (0,0) rectangle (8,8);
\draw[fill=blue!30] (24,24) rectangle (32,32);
\draw (4,4) node {\large 49};
\draw (28,28) node {\large 49};

\end{tikzpicture}
	\caption{\label{fig:nestedVh}Illustration of three $\calH^2$-approximation spaces on $D\times D=[0,1]^2$ for three binary, nested, and perfectly balanced cluster trees. No approximation is performed within the red blocks. The blue blocks are approximated by tensorized iterated interpolation with the inscribed numbers indicating the total degrees of freedom for the specific blocks. $\beta=3$, $\alpha=2$, and $\delta=1$ were used as parameters in \cref{eq:rankdistribution} for this example. The $\calH^2$-approximation spaces are not nested, but have a similar structure which leads to an approximate multi-level hierarchy.}
\end{figure}
To enhance clarity, we write $\calP_t^{\text{pw},\ell}=\calP_t^{\text{pw}}$ and $\calP_{t\times s}^{\text{pw},\ell}=\calP_{t\times s}^{\text{pw}}$ to indicate the level dependence of the polynomial spaces and introduce the operator
\[
\Pi_{h,\ell}^{\calH}=\Pi^{\calH_\ell}\Pi_{h_\ell}^{\mix}.
\]
We denote the range of the operator as
\begin{equation}\label{eq:W}
W^{\calH}_\ell=\Pi^\calH\big(V_h\otimes V_h)\subset L^2(D\times D),\qquad\ell=0,1,\ldots,L.
\end{equation}
In analogy to our previous considerations, these spaces are not nested as well.

\subsubsection{The estimator}
With an approximate multilevel hierarchy of the approximation spaces $W_\ell^\calH$ at hand, our multilevel estimator is defined as follows.
\begin{definition}\label{def:mlmc}
	Setting $\Pi_{h_{-1}}^{\mix}g=0$, we define the \emph{$\calH^2$-formatted multilevel sample covariance estimator ($\calH^2$-MLSCE)} as
	\begin{align}\label{eq:mlmc}\tag{$\calH^2$-MLSCE}
		\bbE[\Pi_{h,L}^{\calH}g]
		\approx
		E_L^{ML}[\Pi_{h,L}^{\calH}g]
		=
		\sum_{\ell=0}^L\Pi^{\calH_L}E_\ell\big[\big(\Pi_{h,\ell}^{\calH}-\Pi_{h,\ell-1}^{\calH}\big)g\big]
	\end{align}
	with the single level estimators
	\[
	E_\ell\big[\big(\Pi_{h,\ell}^{\calH}-\Pi_{h,\ell-1}^{\calH}\big)g\big]
	=
	\frac{1}{M_{\ell}}\sum_{k=1}^{M_{\ell}}\big(\Pi_{h,\ell}^{\calH}-\Pi_{h,\ell-1}^{\calH}\big)\Big(z^{(k)}\otimes z^{(k)}\Big),\quad\ell=0,\ldots,L,
	\]
	given by i.i.d.\ samples $z^{(k)}$, $k=1,\ldots,M_\ell$, $M_\ell\in\bbN$, of the centered Gaussian random field induced by $g$.
\end{definition}
Implementation of this estimator is nontrivial and is discussed in the next section.

\subsubsection{Algorithmic realization}

Careful investigation of the $\calH^2$-MLSCE reveals that
\[
\big(\Pi_{h,\ell}^{\calH}-\Pi_{h,\ell-1}^{\calH}\big)\colon L^2(D\times D)\to W_\ell^\calH+W_{\ell-1}^\calH,
\]
for which an efficient algorithmic treatment is not immediate. A representation in $W_\ell^\calH$ or $W_{\ell-1}^\calH$ would be preferable. To avoid this issue we may rewrite $\calH^2$-MLSCE as
\[
E_L^{ML}[\Pi_{h,L}^{\calH}g]
=
\sum_{\ell=0}^L\Pi^{\calH_L}E_\ell\big[\Pi_{h,\ell}^{\calH}g\big]
-
\sum_{\ell=0}^{L-1}\Pi^{\calH_L}E_{\ell+1}\big[\Pi_{h,\ell}^{\calH}g\big],
\]
while understanding that single level estimators $E_\ell[\cdot]$ with the same index are taken from the same samples. An algorithm implementing the $\calH^2$-MLSCE in this form is given in \cref{alg:H2MLSCE}.

\begin{algorithm}
\caption{\label{alg:H2MLSCE}Algorithmic implementation of the $\calH^2$-MLSCE}
\begin{algorithmic}[1]
\State \textbf{Input:} Finite element space $V_{h_0}$, sample numbers $M_0,\ldots,M_L$, $\alpha$, $\beta$, $\delta$
\State \textbf{Output:} Coefficients $\{\bfU_{t\times s}^{(L)}\}_{t\times s\in L_{I_L\times I_L}}$ of the $\calH^2$-MLSCE $E_L^{ML}[\Pi_{h,L}^{\calH}g]\in W_L^\calH$.
\bigskip
\State Compute cluster tree $T_{I_0}$ and block cluster tree $T_{{I_0}\times {I_0}}$ using standard codes.
\State Generate $V_{h_\ell}$, $\ell=1,\ldots,L$ through uniform refinement.
\State Generate $T_{I_\ell}$ using \cref{alg:nestedct} and $T_{{I_\ell}\times {I_\ell}}$ using standard codes for $\ell=1,\ldots,L$.
\State Compute moment matrices $\{\bfM_r\}_{r\in L_{I_\ell}}$, $\ell=0,\ldots,L$, using standard codes.
\State Compute transfer matrices $\{\bfE_r\}_{r\in T_{I_\ell}}$, $\{\bfF_r\}_{r\in T_{I_\ell}}$, $\ell=0,\ldots,L$, using standard codes.
\State Compute matrices $\{\bfQ_r\}_{r\in L_{I_L}}$, $\ell=0,\ldots,L$.
\smallskip
\For{$\ell=0,\ldots,L$}
\State Compute coefficients $\{\bfU_{t\times s}^{(\ell)}\}_{t\times s\in L_{I_\ell\times I_\ell}}$ of $E_\ell\big[\Pi_{h,\ell}^{\calH}g\big]$ using the $\calH^2$-SCE.
\If{$\ell\neq L$}
\State Compute coefficients $\{\bfU_{t\times s}^{(\ell,1)}\}_{t\times s\in L_{I_\ell\times I_\ell}}$ of $E_{\ell+1}\big[\Pi_{h,\ell}^{\calH}g\big]$ using the $\calH^2$-SCE
\State Update $\{\bfU_{t\times s}^{(\ell)}\}_{t\times s\in L_{I_\ell\times I_\ell}}-\!\!=\{\bfU_{t\times s}^{(\ell,1)}\}_{t\times s\in L_{I_\ell\times I_\ell}}$ using standard codes.
\EndIf
\EndFor
\State Compute coefficients $\{\bfU_{t\times s}^{(L)}\}_{t\times s\in L_{I_L\times I_L}}$ of the $\calH^2$-MLSCE using \cref{alg:MLreduction}.
\end{algorithmic}
\end{algorithm}

The key difference of the algorithm compared to other multilevel Monte Carlo algorithms is the efficient computation of the final estimator from the single level Monte Carlo estimators in line 13. This modification is necessary, since the approximation spaces $W_\ell^\calH$ are not nested, cf.\ our discussion at the end of \cref{sec:mlhier}. The efficient computation of the final estimator is performed in \cref{alg:MLreduction}. We derive this algorithm in \cref{sec:MLMCalgo} and show in \cref{lem:reductioncomplexity} that it has complexity
\begin{equation}\label{eq:MLreductioncompintro}
\calO\bigg(\frac{(\alpha+\beta)^{\lceil 3\delta d\rceil}}{(\alpha+\beta)^{\delta d}}|I_L|\bigg).
\end{equation}
The matrices $\bfJ_{t}$ in the algorithm are to be understood as the matrix representation of $\calJ_\ell|_{D_t}\colon V_{h_\ell}|_{D_t}\to V_{h_{\ell+1}}|_{D_t}$, where $\calJ_\ell\colon V_{h_\ell}\to V_{h_{\ell+1}}$ is the canonical prolongation operator between nested finite element spaces.

\begin{algorithm}
\caption{\label{alg:MLreduction}The multilevel $\calH^2$-reduction algorithm}
\begin{algorithmic}[1]
\State \textbf{Input:} Coefficients $\{\bfU_{t\times s}^{(\ell)}\}_{L_{I_\ell\times I_\ell}}$, $\ell=0,\ldots,L$, computed in \cref{alg:H2MLSCE}.
\State \textbf{Output:} Coefficients $\{\bfU_{t\times s}^{(L)}\}_{t\times s\in L_{I\times I}}$ of $E_L^{ML}[\Pi_{h,L}^{\calH}g]\in W_L^\calH$.
\bigskip
\State Set $\bfR_t=\bfQ_t$, $\bfN_t=\bfM_t$ for all $t\in L_{I_L}$.
\State Compute $\bfQ_r=\sum_{r'\in\children(r)}\bfE_{r'}^\intercal \bfQ_{r'}\bfE_{r'}$ for all $r\in T_{I_L}\setminus L_{I_L}$ recursively.
\State Compute $\bfR_t=\sum_{t'\in\children(t)}\bfE_{t',L}^\intercal \bfR_{t'}\bfF_{t'}$ for all $t\in T_{I_L}\setminus L_{I_L}$ recursively.
\State Compute $\bfN_t=\sum_{t'\in\children(t)}\bfE_{t',L}^\intercal \bfN_{t'}\bfJ_{t'}^\intercal$ for all $t\in T_{I_L}\setminus L_{I_L}$ recursively.
\For{$\ell=0,\ldots,L-1$}
\State Set $\bfR_t^{(L,\ell)}=\bfR_t$, $\bfN_t^{(L,\ell)}=\bfN_t$ for all $t\in L_{I_\ell}$,
\State Compute $\bfR_t^{(L,\ell)}=\sum_{t'\in\children(t)}\bfE_{t',L}^\intercal \bfR_{t'}^{(L,\ell)}\bfE_{t',\ell}$ for all $t\in T_{I_\ell}\setminus L_{I_\ell}$ recursively.
\State Compute $\bfN_t^{(L,\ell)}=\sum_{t'\in\children(t)}\bfE_{t',L}^\intercal \bfN_{t'}^{(L,\ell)}\bfJ_{t'}$ for all $t\in T_{I_\ell}\setminus L_{I_\ell}$ recursively.
\State Compute prolongation matrices $\{\bfJ_t^{(\ell)}\}_{t\in T_{I_\ell}}$ using finite element codes.
\EndFor
\bigskip
\For{$\ell=0,\ldots,L-1$}
\smallskip
\For{$t\times s\in L_{I_\ell\times I_\ell}$}
\If{$t\times s\in L_{I_L\times I_L}^+$}
\State Solve
\[
\bfQ_{t}\bfu_{t\times s}^{(L)}\bfQ_{s}^\intercal
=
\begin{cases}
\bfR_{t}^{(L,\ell)}\bfU_{t\times s}^{(\ell)}\Big(\bfR_{s}^{(L,\ell)}\Big)^\intercal& \text{if}~t\times s\in L_{I_\ell\times I_\ell}^+,\\
\bfN_{t}^{(L,\ell)}\bfU_{t\times s}^{(\ell)}\Big(\bfN_{s}^{(L,\ell)}\Big)^\intercal& \text{if}~t\times s\in L_{I_\ell\times I_\ell}^-.\\
\end{cases}
\]
\State Set $\bfU_{t\times s}^{(L)}+\!\!=\bfu_{t\times s}^{(L)}$.
\Else
\State Prolongate $\widetilde{\bfU}_{t\times s}^{(\ell)}=\bfJ_{t}^\intercal\bfU_{t\times s}^{(\ell)}|_{D_t\times D_s}\bfJ_{s}$
\For{$t'\times s'\in\children(t\times s)$}
\If{$t'\times s'\in L_{I_{\ell+1}\times I_{\ell+1}}^+$}
\State Solve $\bfQ_{t'}\bfu_{t'\times s'}^{(L)}\bfQ_{s'}^\intercal=\bfN_{t'}^{(L,\ell)}\widetilde{\bfU}_{t'\times s'}^{(\ell)}\Big(\bfN_{s'}^{(L,\ell)}\Big)^\intercal$
\State Set $\bfU_{t'\times s'}^{(L)}+\!\!=\bfu_{t'\times s'}^{(L)}$
\Else
\State Set $\bfU_{t'\times s'}^{(\ell+1)}+\!\!=\widetilde{\bfU}_{t'\times s'}^{(\ell)}$
\EndIf
\EndFor
\EndIf
\EndFor
\EndFor
\end{algorithmic}
\end{algorithm}

\subsubsection{Approximation result}
We obtain the following approximation result for the $\calH^2$-MLSCE.
\begin{theorem}\label{thm:MLestimatorerror}
	Let $2\theta<0$ and let $g\in H^\gamma(D\times D)$ be $G^\delta(C_G,A,2\theta)$-asymptotically smooth. Then there are $\tilde{\rho}\in(0,1)$ (specified in \cref{eq:rhotilde}), $\alpha\in\bbN$, and $\beta_0\in\bbN$ such that $\zeta^{-2\theta}\tilde{\rho}^{\alpha}<1$ and
	\begin{align*}
		\big\|g-E_L^{ML}[\Pi_{h,L}^{\calH}g]\big\|_{L_{\bbP}^2(\Omega;L^2(D\times D))}
		&\leq
		\frac{C_{\text{lc}}C_{\text{sp}}\tilde{\rho}^\beta}{1-\zeta^{-2\theta}\tilde{\rho}^\alpha}\bigg(h_{\calH,L}^{-2\theta}+(1+2^{-2\theta})\sum_{\ell=0}^L\frac{h_{\calH,\ell}^{-2\theta}}{\sqrt{M_\ell}}\bigg)\\
		&
		\qquad+C_{L^2}^\otimes\bigg( h_L^{{\gamma}}+ (1+2^{\gamma})\sum_{\ell=0}^L\frac{h_\ell^{{\gamma}}}{\sqrt{M_\ell}}\bigg)\|g\|_{H^\gamma(D\times D)}.
	\end{align*}
\end{theorem}
\begin{proof}
	The estimate is proven in the usual way, see, e.g., \cite{BSZ2011}, and we omit the proof for brevity.
\end{proof}

Thus, by increasing $\beta$, $L$, and $M_\ell$, we can approximate the covariance function $g$ to arbitrary precision with the $\calH^2$-SCE. For brevity, we may shorten the theorem to the following.

\begin{corollary}\label{cor:shortMLestimatorerror}
	Let the assumptions of \cref{thm:MLestimatorerror} hold and let $\tilde{\gamma}=\min\{-2\theta,\gamma\}$. Then there are $\tilde{\rho}\in(0,1)$ (specified in \cref{eq:rhotilde}), $\alpha\in\bbN$, $\beta_0\in\bbN$, and a constant
	\[
	0<C_{\text{MLE}}=C_{\text{MLE}}\Big(C_{\text{lc}}C_{\text{sp}}\tilde{\rho}^{\beta_0},\zeta^{-2\theta}\tilde{\rho}^\alpha,C_{L^2}^\otimes,-2\theta,\gamma,\|g\|_{H^\gamma(D\times D)}\Big)
	\]
	such that $\zeta^{-2\theta}\tilde{\rho}^{\alpha}<1$ and
	\[
	\big\|g-E_L^{ML}[\Pi_{h,L}^{\calH}g]\big\|_{L_{\bbP}^2(\Omega;L^2(D\times D))}
	\leq
	C_{\text{MLE}}\bigg(h_L^{\tilde{\gamma}}+\sum_{\ell=0}^L\frac{h_\ell^{\tilde{\gamma}}}{\sqrt{M_\ell}}\bigg)
	\]
	for all $\beta\geq\beta_0$.
\end{corollary}
\begin{proof}
	The specific construction of $\{\mathcal{T}_{h_\ell}\}_{\ell=0}^L$, $\{T_{I_\ell}\}_{\ell=0}^L$, and $\{V_{h_\ell}\}_{\ell=0}^L$ using uniform refinement implies that $c^{-1}h_\ell\leq h_{\calH,\ell}\leq ch_\ell$ for $\ell=0,\ldots, L$. This yields the assertion.
\end{proof}

In analogy to \cref{lem:fembilinearerror} we obtain the following bound on the bilinear form induced by the covariance function.

\begin{corollary}\label{thm:blfestimate}
	Under the assumptions of \cref{cor:shortMLestimatorerror} there is $\beta_0\in\bbN$ such that
	\begin{align*}
		&\bigg\|\int_{D}\int_{D}\Big(g(\bfx,\bfy)-\Pi^{\calH_L}E_L^{ML}[\Pi_{h_L}^{\mix}g(\bfx,\bfy)]\Big)u_h(\bfx)v_h(\bfy)\dd\mu(\bfx)\dd\mu(\bfy)\bigg\|_{L_{\bbP}^2(\Omega)}\\
		&\qquad\qquad\qquad\qquad\qquad\qquad\leq C_{\text{MLE}}\bigg(h_L^{\tilde{\gamma}}+\sum_{\ell=0}^L\frac{h_\ell^{\tilde{\gamma}}}{\sqrt{M_\ell}}\bigg)\|u_h\|_{L^2(D)}\|v_h\|_{L^2(D)},
	\end{align*}
	for all $\beta\geq\beta_0$ with $\tilde{\rho}$ as in \cref{eq:rhotilde}.
\end{corollary}

\subsubsection{Choice of sample numbers}

Our final remarks are concerned with choosing the sample numbers of the $\calH^2$-MLSCE.
Combining \cref{eq:H2SCEcompl} and \cref{eq:MLreductioncompintro} yields that the $\calH^2$-MLSCE can be computed in
\[
\calO\bigg((\alpha+\beta)^{2\delta d}\sum_{\ell=0}^LM_\ell|I_\ell|+\frac{(\alpha+\beta)^{\lceil 3\delta d\rceil}}{(\alpha+\beta)^{\delta d}}|I_L|\bigg)
\]
operations. Thus, it remains to choose the sample numbers such that accuracy of the finest level is achieved with minimal work. In complete analogy to various references, we mention \cite[Appendix D]{HHKS2021} or \cite{KSS2015} for example, we state the following theorem without proof.
\begin{theorem}\label{thm:costcomplexity}
	Let the assumptions of \cref{cor:shortMLestimatorerror} hold and choose $\varepsilon>0$. The $\calH^2$-MLSCE with
	\[
	L=\frac{d}{\tilde{\gamma}}\bigg|\frac{\log(\varepsilon^{-1})}{\log(C_{\text{uni}})}\bigg|
	\]
	and sample numbers
	\[
	M_\ell=M_0C_{\text{uni}}^{-2\ell (1+\tilde{\gamma}/d)/3},\qquad\ell=0,\ldots,L,
	\]
	with
	\[
	M_0
	=
	\begin{cases}
		C_{\text{uni}}^{2\tilde{\gamma}L/d}&\text{for}~2\tilde{\gamma}>d,\\
		C_{\text{uni}}^{2\tilde{\gamma}L/d}L^2&\text{for}~2\tilde{\gamma}=d,\\
		C_{\text{uni}}^{2(1+\tilde{\gamma}/d)L/3}&\text{for}~2\tilde{\gamma}<d,
	\end{cases}
	\]
	achieves error estimates
	\[
	\big\|g-E_L^{ML}[\Pi_{h,L}^{\calH}g]\big\|_{L_{\bbP}^2(\Omega;L^2(D\times D))}=\calO(\varepsilon)
	\]
	and
	\[
	\sup_{u_h,v_h\in V_{h_L}}
	\frac{\Big\|\int_{D}\int_{D}\big(g(\bfx,\bfy)-\Pi^{\calH_L}E_L^{ML}[\Pi_{h_L}^{\mix}g(\bfx,\bfy)]\big)u_h(\bfx)v_h(\bfy)\dd\mu(\bfx)\dd\mu(\bfy)\Big\|_{L_{\bbP}^2(\Omega)}}{\|u_h\|_{L^2(D)}\|v_h\|_{L^2(D)}}=\calO(\varepsilon)
	\]
	in a computational complexity of
	\[
	\begin{cases}
		\calO(\varepsilon^{-2})&\text{for}~2\tilde{\gamma}>d,\\
		\calO\big(\varepsilon^{-2}|\log(\varepsilon^{-1})|^3\big)&\text{for}~2\tilde{\gamma}=d,\\
		\calO(\varepsilon^{-d/\tilde{\gamma}})&\text{for}~2\tilde{\gamma}<d.\\
	\end{cases}
	\]
\end{theorem}
Thus, for $2\tilde{\gamma}>d$, the overall error is dominated by the Monte Carlo error, whereas for $2\tilde{\gamma}<d$ the overall error is dominated by the error of the approximation spaces $V_{h_l}$.

We note that these computational complexities are in line with the wavelet-based approach from \cite{HHKS2021}, but the $\calH^2$-approach does not require a hierarchical basis. In contrast, wavelet-based approaches are theoretically also applicable if the smoothness of the kernel function is finite, which is, see also \cref{rem:finitesmoothness}, asymptotically not the case for the $\calH^2$-approach due to the increasingly higher polynomial degrees required for interpolation.

\section{Proofs of $\calH^2$-approximation properties}\label{sec:Gevreyapprox}

\subsection{Proof of \cref{thm:abGevreyInterpolationintro}}

The main ingredient to show \cref{thm:abGevreyInterpolationintro} is the following.
\begin{theorem}\label{thm:m11GevreyInterpolation}
Let $f\in G^\delta([-1,1],C_G,A)$, $\rho(r)=r+\sqrt{1+r^2}$, and $m\in\bbN$, $m\geq 3$. Then it holds
\[
\min_{q\in\calP_m}\|f-q\|_{C([-1,1])}\leq C(A,\delta)C_G\rho(1/A)^{-m^{1/\delta}/e^2},
\]
where $C(A,\delta)$ is monotonically increasing in $A$.
\end{theorem}
\begin{proof}
The proof is inspired by the one of \cite[Proposition 4.1]{OSZ2022}. Denote by $I_3\colon H^2([-1,1])\to\calP_3$ the Hermite interpolation operator given by $I_3f(\pm1)=f(\pm1)$, $(I_3f)'(\pm1)=f'(\pm1)$ and, for $m\in\bbN$, $m\geq 3$, denote by $\pi_{m-2,0}\colon L^2([-1,1])\to\calP_{m-2}$ the $L^2$-orthogonal projection onto the first $m-1$ Legendre polynomials. Then, the projector $H^2([-1,1])\to\calP_m$ defined by
\[
(\pi_{m,2}f)(x)=(I_3f)(x)+\int_{-1}^{x}\int_{-1}^y\big(\pi_{m-2,0}\big((f-I_3f)''\big)\big)(z)\dd z\dd y
\]
satisfies the error estimate, see \cite[Theorem A.1]{CDS2005},
\[
\|f-\pi_{m,2}f\|_{H^2([-1,1])}^2\leq C\frac{(m-1-k)!}{(m-1+k)!}\big\|f^{(k+2)}\big\|_{L^2([-1,1])}^2,\qquad 2\leq k\leq m-1.
\]

Now, fix $\alpha=(2\rho(1/A)^{\rho(A)} A)^{-1/\delta}$, $k=\lfloor\alpha \gamma m^{1/\delta}\rfloor$ with $\gamma=\min\{\max\{\frac{2}{\alpha m^{1/\delta}},1\},\frac{m-1}{\alpha m^{1/\delta}}\}$, and note that $2\leq k\leq m-1$, $k\leq\alpha\gamma m^{1/\delta}\leq k+1$, and $\rho(1/A)^{\rho(A)}\leq(2/A+1)^{2A+1}\defis\Xi(A)$. Gevrey regularity $f\in G^\delta([-1,1],C_G,A)$ and Stirling's formula $\sqrt{2\pi n}(n/e)^n\leq n!\leq e\sqrt{n}(n/e)^n$, $n\in\bbN$ imply
\begin{align*}
\|f-\pi_{m,2}f\|_{H^2([-1,1])}^2
\leq{}& CC_G^2A^{2k+4}\frac{(m-1-k)!}{(m-1+k)!}\big((k+2)!\big)^{2\delta}\\
\leq{}& CC_G^2A^{2k+4}\frac{e^{1+2k}}{\sqrt{2\pi}}\frac{(m-1-k)^{m-1-k+1/2}}{(m-1+k)^{m-1+k+1/2}}\big(k!(k+2)^2\big)^{2\delta}\\
\leq{}& CC_G^2A^{2k+4}\frac{e^{1+2k}}{\sqrt{2\pi}}\frac{(m-1-k)^{m-1-k+1/2}}{(m-1+k)^{m-1+k+1/2}}e^{2\delta(1-k)}k^{2k\delta}k^\delta(k+2)^{4\delta}\\
\leq{}& CC_G^2A^{2k+4}\frac{e^{1+2\delta+2(1-\delta)k}}{\sqrt{2\pi}}\bigg(\frac{m-1-k}{m-1+k}\bigg)^{m-1-k+1/2}m^{-2k}k^{2k\delta}k^\delta(k+1)^{4\delta}.
\end{align*}
Since $1-\delta\leq 0$, $m-1-k+1/2\geq 0$ for $k\leq m-1$, $m^{-k}\leq (\alpha\gamma/k)^{\delta k}$, and $k^{\delta}(k+2)^{4\delta}\leq C(\delta)2^{2k}$ for $k\geq 2$ this implies
\[
\|f-\pi_{m,2}f\|_{H^2([-1,1])}\leq C(\delta)C_GA^{k+2}\gamma^{\delta k}\alpha^{\delta k}2^k.
\]
We next remark that $\gamma^{\delta k}\leq 1$ for $2\leq\alpha m^{1/\delta}$. For $2>\alpha m^{1/\delta}$, we remark that $\gamma^{\delta k}\leq\gamma^{2\delta}\leq C(\delta)(\Xi(A)A)^2$, where $\Xi(A)A$ is continuous and monotonically increasing on $(0,\infty)$ with $\lim_{t\to 0}\Xi(t)t=2$. Thus, $\gamma^{\delta k}\leq\chi(A,\delta)$ is monotonically increasing in $A$ with $\chi(A,\delta)\geq 4C(\delta)$. The continuous embedding $H^2([-1,1])\hookrightarrow L^\infty([-1,1])$ and the definition of $\alpha$ then yield
\[
\|f-\pi_{m,2}f\|_{C([-1,1])}
\leq
C(A,\delta)C_GA^2\rho(1/A)^{\rho(A)}\rho(1/A)^{-\rho(A)(k+1)}
\leq
C(A,\delta)C_G\rho(1/A)^{-\rho(A)\alpha\gamma m^{1/\delta}},
\]
where $C(A,\delta)$ is monotonically increasing in $A$. To obtain the desired exponent, we consider that $\rho(A)\alpha(A,\delta)$ is monotonically increasing in $\delta$, and that it is bounded from below by $e^{-2}$ for $\delta = 1$. For $\alpha m^{1/\delta}< m-1$ this yields $\rho(A)\alpha\gamma m^{1/\delta}\geq m^{1/\delta}/e^2$ due to $\gamma\geq 1$. For $\alpha m^{1/\delta}\geq m-1$ we observe that $\rho(A)\alpha\gamma m^{1/\delta}\geq m-1\geq m^{1/\delta}/e^2$, which yields the assertion.
\end{proof}
The following corollary generalizes the result to arbitrary intervals.
\begin{corollary}\label{cor:abGevreyInterpolation}
For any $f\in G^\delta([a,b],C_G,A)$, $B=A(b-a)/2$, and $m\in\bbN$, $m\geq 3$, it holds that
\[
\min_{q\in\calP_m}\|f-q\|_{C([a,b])}\leq C(B,\delta)C_G\rho(1/B)^{-m^{1/\delta}/e^2},
\]
where $C(B,\delta)$ is monotonically increasing in $B$.
\end{corollary}
\begin{proof}
Denoting $\Phi^{[a,b]}\colon [-1,1]\to[a,b]$ with $\Phi^{[a,b]}(t)=(b+a)/2+t(b-a)/2$, one easily verifies that $f\in G^\delta([a,b],C_G,A)$ implies $f\circ\Phi^{[a,b]}\in G^\delta([-1,1],C_G,B)$. \cref{thm:m11GevreyInterpolation} yield the assertion.
\end{proof}
Thus, we have shown \cref{thm:abGevreyInterpolationintro}.

\subsection{Interpolation of Gevrey functions}
The remainder of this section aims at proving \cref{thm:kernell2errorintro}. To this end, we consider polynomial interpolation $\calI_m^{[a,b]}\colon C([a,b])\to \calP_m$ on $m+1$ distinct points in $[a,b]$. To avoid technicalities we assume that the interpolation points are Chebyshev points. This implies that the interpolation operator is stable, i.e.,
\[
\big\|\calI_m^{[a,b]}[f]\big\|_{C([a,b])}
\leq
\Lambda_m \|f\|_{C([a,b])},
\]
for all $m\in\bbN$, with stability constant $1\leq\Lambda_m\leq\frac{2}{\pi}\ln(m)+1$, see, e.g., \cite[Theorem 1.2]{Riv1974}. Moreover, it is well known that 
for $m\in\bbN$ and $f\in C([a,b])$ it holds
\begin{equation}\label{eq:interpolationoptimal}
\big\|f-\calI_m^{[a,b]}[f]\big\|_{C([a,b])}
\leq
(\Lambda_m+1)\min_{q\in\calP_m}\|f-q\|_{C([a,b])},
\end{equation}
i.e., that we can bound the polynomial interpolation error in terms of the best approximation error among all polynomials with the same degree.

For $Q=\bigtimes_{i=1}^d[a_i,b_i]$, $f\in C(Q)$, and $m\in\bbN$, we define the tensor product interpolation operator on $Q$ as $\calI_m^Q=\bigotimes_{i=1}^d\calI_{m,i}^{[a_i,b_i]}$. Combining \cref{eq:interpolationoptimal} with \cref{cor:abGevreyInterpolation} we obtain the following approximation result for interpolation on tensor product domains.
\begin{theorem}\label{thm:QGevreyInterpolation}
Let $Q=\bigtimes_{i=1}^d[a_i,b_i]$, $f\in G^\delta(Q,C_G,A)$, and $m\in\bbN$, $m\geq 3$. Then it holds
\[
\|f-\calI_m^Q[f]\|_{C(Q)}\leq  C(A\diam_\infty(Q)/2,\delta) C_Gd(\Lambda_m+1)^d\rho\bigg(\frac{2}{A\diam_\infty(Q)}\bigg)^{-m^{1/\delta}/e^2},
\]
where $C(A,\delta)$ is monotonically increasing in $A$.
\end{theorem}
\begin{proof}
In complete analogy to the proof of \cite[Corollary 4.21]{Bor2010}, using \cref{cor:abGevreyInterpolation}.
\end{proof}

\subsection{Interpolation of Gevrey kernels}
As outlined in \cref{sec:Gdelta} we are interested in the approximation of functions satisfying the regularity estimate \cref{eq:introgevreysmoothkernel}.
Using the approximation estimate from \cref{thm:QGevreyInterpolation}, the following theorem generalizes the very similar result for asymptotically smooth kernels proven in \cite[Theorem 4.22]{Bor2010}.
\begin{theorem}\label{thm:admissibleapproximation}
Let $Q_t=\bigtimes_{i=1}^d[a_i,b_i]$ and $Q_s=\bigtimes_{i=d+1}^{2d}[a_i,b_i]$. Let $\eta>0$ and $Q_t$ and $Q_s$ be \emph{admissible}, i.e.,
Let $g$ be $G^{\delta}(C_G,A)$-asymptotically smooth on $Q_t\times Q_s$ and $\tilde{g}=\calI_m^{Q_t\times Q_s}[g]$. Then it holds for $m\in\bbN$, $m\geq 3$,
\begin{align}\label{eq:gevreycontraction}
\|g-\tilde{g}\|_{C(Q_t\times Q_s)}\leq C(A\eta,\delta)C_G\frac{2d(\Lambda_m+1)^{2d}}{\dist_2(Q_t,Q_s)^{2\theta+d}}\rho\bigg(\frac{1}{A\eta}\bigg)^{-m^{1/\delta}/e^2}.
\end{align}
\end{theorem}
\begin{proof}
In complete analogy to the proof of \cite[Theorem 4.22]{Bor2010}, using \cref{thm:QGevreyInterpolation} rather than approximation estimates for analytic functions.
\end{proof}
To improve readability we may note that $\lim_{t\to\infty}q(t)\tilde{\rho}^{t^{1/\delta}}=0$ for any polynomial $q$ and $\tilde{\rho}\in(0,1)$ to follow \cite[Remark 4.23]{Bor2010} and reformulate \cref{eq:gevreycontraction} in \cref{thm:admissibleapproximation} as
\begin{align}\label{eq:rhotilde}
\|g-\tilde{g}\|_{C(Q_t\times Q_s)}
\leq
\frac{C_{\text{in}}}{\dist_2(Q_t,Q_s)^{2\theta+d}}\tilde{\rho}^{m^{1/\delta}},
\qquad
\tilde{\rho}\isdef\min\bigg\{\frac{A\eta}{A\eta+1},\frac{A\eta}{2}\bigg\}^{1/e^2}>\rho\bigg(\frac{1}{A\eta}\bigg)^{-1/e^2},
\end{align}
for some fixed $C_{\text{in}}>0$. Here, we have used implicitly that $\Lambda_m\leq m$.

\begin{remark}\label{rem:finitesmoothness}
The classical results for asymptotically smooth kernel functions depend on the analyticity of the kernel function in admissible clusters since these estimates are based on analytic continuations into Bernstein ellipses in the complex plane. In contrast, the arguments of our generalizations to $G^{\delta}(C_G,A)$-asymptotically smooth kernels only require finite smoothness in each direction and do not require extensions into the complex plane.
\end{remark}

All further results from the variable-order $\calH^2$-theory for asymptotically smooth kernels can be generalized with minor modifications. In the following we use the common assumptions and state a slightly modified error estimate in the $L^2$-norm, rather than the maximums norm.

\subsection{$L^2$-error of variable-order $\calH^2$-approximations}
For Gevrey-regular kernels, the approximation error in each admissible block-cluster can be estimated as follows.
\begin{corollary}
\label{cor:blockiteratedinterpolationestimate}
Let $2\theta<0$ and let the kernel function $g\colon\bbR^d\times\bbR^d\to\bbR$ be $G^\delta(C_G,A,2\theta)$-asymptotically smooth, and let $\alpha\in\bbN_0$. Then there are constants $C_{\text{in}}\in\bbR_{>0}$ and $\beta_0\in\bbN_0$ such that
\[
\big\|g-\calI_{t_0\times s_0}^{t\times s}g\big\|_{C(Q_{t_0}\times Q_{s_0})}
\leq
C_{\text{in}}
\bigg(\frac{\tilde{\rho}^{\beta+\alpha(p-\level(t))}}{\diam_\infty(Q_t)^{2\theta+d}}\bigg)^{1/2}
\bigg(\frac{\tilde{\rho}^{\beta+\alpha(p-\level(s))}}{\diam_\infty(Q_s)^{2\theta+d}}\bigg)^{1/2}
\]
holds with $\tilde{\rho}$ as in \cref{eq:rhotilde} for all $\beta\geq\beta_0$, all blocks $t\times s\in L_{I\times I}^+$ satisfying \cref{eq:admissibility}, and all $t_0\in L_t$, $s_0\in L_s$.
\end{corollary}
\begin{proof}
The proof follows the arguments of \cite{BLM2005} and \cite[Chapter 4.7]{Bor2010} with only minor modifications.
\end{proof}
At first sight, the error estimate of \cref{cor:blockiteratedinterpolationestimate} (and the ones from the following results) seems to be independent of the Gevrey class $\delta\geq 1$. This is true for the upper bound. However, we note that the polynomial degrees chosen as in \cref{eq:kidelta} are $\delta$-dependent, implying that the interpolation operator $\calI_{t_0\times s_0}^{t\times s}$ is $\delta$-dependent as well. The choice to modify the polynomial degrees, rather than the error estimates, will simplify the analysis later on.

Although the results from the literature can be generalized to Gevrey kernels, most of the analysis in the literature is based on a $C(Q_{t_0}\times Q_{s_0})$-type estimate, which is not compatible with the $L^2$-setting of the Monte Carlo type error analysis, for which an $L^2$-estimate is preferable.

\begin{corollary}
\label{cor:L2blockiteratedinterpolationestimate}
Let the assumptions of \cref{sec:approxGvoH2} and the assumptions of \cref{cor:blockiteratedinterpolationestimate} hold. Then it holds
\begin{align*}
\bigg\|g-\sum_{t_0\times s_0\in L_{t\times s}}\calI_{t_0\times s_0}^{t\times s}g\bigg\|_{L^2(D_t\times D_s)}
\leq
C_{\text{lc}} h_{\calH}^{-2\theta}\tilde{\rho}^\beta(\zeta^{-2\theta}\tilde{\rho}^\alpha)^{p-\level(t)/2-\level(s)/2},
\end{align*}
where $C_{\text{lc}}=C_{\text{in}}C_{\text{cu}}C_{\text{gr}}^{-2\theta}$.
\end{corollary}
\begin{proof}
The assumptions of \cref{sec:approxGvoH2} imply
\[
\diam_{\infty}(Q_t)\leq\zeta^{\level(t')-\level(t)}\diam_{\infty}(Q_{t'})\leq C_{\text{gr}}h_{\calH}\zeta^{p-\level(t)}
\]
for all $t'\in L_t$, $t\in T_I$. Thus,
\[
\frac{\mu(D_t)}{\diam_{\infty}(Q_t)^{2\theta+d}}
\leq
\frac{C_{\text{cu}}}{\diam_{\infty}(Q_t)^{2\theta}}
\leq
C_{\text{cu}}C_{\text{gr}}^{-2\theta}h_{\calH}^{-2\theta}\zeta^{-2\theta(p-\level(t) )}.
\]
The assertion follows from H\"older's inequality and \cref{cor:blockiteratedinterpolationestimate} due to
\begin{align*}
\bigg\|g-\sum_{t_0\times s_0\in L_{t\times s}}\calI_{t_0\times s_0}^{t\times s}g\bigg\|_{L^2(D_t\times D_s)}
&\leq
\max_{t_0\times s_0\in L_{t\times s}}\|g-\calI_{t_0\times s_0}^{t\times s}g\|_{C(Q_{t_0}\times Q_{s_0})}\mu(D_t)^{1/2}\mu(D_s)^{1/2}\\
&\leq
C_{\text{in}}
\bigg(\frac{\mu(D_t)\tilde{\rho}^{\beta+\alpha(p-\level(t))}}{\diam_\infty(Q_t)^{2\theta+d}}\bigg)^{1/2}
\bigg(\frac{\mu(D_s)\tilde{\rho}^{\beta+\alpha(p-\level(s))}}{\diam_\infty(Q_s)^{2\theta+d}}\bigg)^{1/2}.
\end{align*}
\end{proof}

\subsection{Proof of \cref{thm:kernell2errorintro}}\label{sec:L2toH2}
The purpose of this subsection is to combine the $L^2$ approximation estimates for block-clusters from \cref{cor:L2blockiteratedinterpolationestimate} to a global estimate in $L^2(D\times D)$. To this end, since the intersection of two block-clusters has measure zero, computing $\Pi^\calH k$ is equivalent to computing the $L^2(D_t\times D_s)$ projections $\Pi_{t\times s}^\calH k$ of $k|_{D_t\times D_s}$ onto $\calP_{t\times s}^{\text{\text{pw}}}$ and setting
\begin{equation}\label{eq:L2local}
\Pi^\calH k=\bigg(\sum_{t\times s\in L_{I\times I}^+}\Pi_{t\times s}^\calH k+\sum_{t\times s\in L_{I\times I}^-}k|_{D_t\times D_s}\bigg)\in V^\calH
\end{equation}
for all $k\in L^2(D\times D)$. Here, we extend $\Pi_{t\times s}^\calH k$ and $k|_{D_t\times D_s}$ by zero outside $t\times s$ to simplify notation. \Cref{cor:L2blockiteratedinterpolationestimate} immediately implies 
\begin{equation}\label{lem:L2blockiteratedinterpolationestimate}
\big\|g-\Pi^\calH g\big\|_{L^2(D_t\times D_s)}
=
\big\|g-\Pi_{t\times s}^\calH g\big\|_{L^2(D_t\times D_s)}
\leq
C_{\text{lc}} h_{\calH}^{-2\theta}\tilde{\rho}^\beta(\zeta^{-2\theta}\tilde{\rho}^\alpha)^{p-\level(t)/2-\level(s)/2}
\end{equation}
for all blocks $t\times s\in L_{I\times I}^+$. For the following approximation estimate we further need the sparsity constant of a block-cluster tree.
Given a block-cluster tree $T_{I\times I}$, its \emph{sparsity constant $C_{\text{sp}}$} is defined as
\[
C_{\text{sp}}=\max_{t\in T_I}\big|\big\{s\in T_I\colon t\times s\in T_{I\times I}\big\}\big|,
\]
cf.\ also \cite{Hac2015}.

\begin{theorem}\label{thm:L2interpolationestimate}
Let the assumptions of \cref{cor:L2blockiteratedinterpolationestimate} hold. Choose $\alpha\in\bbN$ such that $\zeta^{-2\theta}\tilde{\rho}^\alpha<1$. Then there is $\beta_0\in\bbN$ such that
\[
\big\|g-\Pi^\calH g\big\|_{L^2(D\times D)}
\leq
\frac{C_{\text{lc}}C_{\text{sp}}h_{\calH}^{-2\theta}\tilde{\rho}^{\beta}}{1-\zeta^{-2\theta}\tilde{\rho}^{\alpha}}
\]
for all $\beta\geq\beta_0$ with $\tilde{\rho}$ as in \cref{eq:rhotilde}.
\end{theorem}
\begin{proof}
Due to $L^2(D\times D)\simeq L^2(D)\otimes L^2(D)\simeq L^2(D;L^2(D))$ we may write
\[
\big\|g-\Pi^{\calH}g\big\|_{L^2(D\times D)}
=
\sup_{\substack{u,v\in L^2(D)\\u,v\neq 0}}
\frac{\int_D\int_D \big(g(\bfx,\bfy)-\Pi^{\calH}g(\bfx,\bfy)\big)u(\bfx)v(\bfy)\dd\mu(\bfx)\dd\mu(\bfy)}{\|u\|_{L^2(D)}\|v\|_{L^2(D)}}.
\]
Using \cref{lem:L2blockiteratedinterpolationestimate}, the Cauchy-Schwartz inequality, and sparsity of $T_{I\times I}$ the numerator is estimated by
\begin{align*}
&\int_D\int_D \big(g(\bfx,\bfy)-\Pi^{\calH}g(\bfx,\bfy)\big)u(\bfx)v(\bfy)\dd\mu(\bfx)\dd\mu(\bfy)\\
&\qquad\leq\sum_{t\times s\in L^+_{I\times I}}\big\|g-\Pi^{\calH}g\big\|_{L^2(D_t\times D_s)}\|u\|_{L^2(D_t)}\|v\|_{L^2(D_s)}\\
&\qquad\leq C_{\text{lc}} h_{\calH}^{-2\theta}\tilde{\rho}^\beta\sum_{t\times s\in L_{I\times I}^+}(\zeta^{-2\theta}\tilde{\rho}^\alpha)^{(p-\level(t))/2}\|u\|_{L^2(D_t)}(\zeta^{-2\theta}\tilde{\rho}^\alpha)^{(p-\level(s))/2}\|v\|_{L^2(D_s)}\\
&\qquad\leq C_{\text{lc}}C_{\text{sp}}h_{\calH}^{-2\theta}\tilde{\rho}^\beta\bigg(\sum_{t\in T_I}(\zeta^{-2\theta}\tilde{\rho}^\alpha)^{p-\level(t)}\|u\|_{L^2(D_t)}^2\bigg)^{1/2}\bigg(\sum_{s\in T_I}(\zeta^{-2\theta}\tilde{\rho}^\alpha)^{p-\level(s)}\|v\|_{L^2(D_s)}^2\bigg)^{1/2}\\
&\qquad\leq C_{\text{lc}}C_{\text{sp}}h_{\calH}^{-2\theta}\tilde{\rho}^\beta\bigg(\sum_{\ell=0}^p(\zeta^{-2\theta}\tilde{\rho}^\alpha)^{p-\ell}\sum_{\substack{t\in T_I\\\level(t)=\ell}}\|u\|_{L^2(D_t)}^2\bigg)^{1/2}\bigg(\sum_{\ell=0}^p(\zeta^{-2\theta}\tilde{\rho}^\alpha)^{p-\ell}\sum_{\substack{s\in T_I\\\level(s)=\ell}}\|v\|_{L^2(D_s)}^2\bigg)^{1/2}.
\end{align*}
Finally, $\zeta^{-2\theta}\tilde{\rho}^\alpha<1$ implies
\[
\sum_{\ell=0}^p(\zeta^{-2\theta}\tilde{\rho}^\alpha)^{p-\ell}\sum_{\substack{t\in T_I\\\level(t)=\ell}}\|u\|_{L^2(D_t)}^2
\leq
\|u\|_{L^2(D)}^2\sum_{\ell=0}^p(\zeta^{-2\theta}\tilde{\rho}^\alpha)^{\ell}
\leq
\frac{\|u\|_{L^2(D)}^2}{1-\zeta^{-2\theta}\tilde{\rho}^\alpha},
\]
which yields the assertion.
\end{proof}
Thus, we have shown \cref{thm:kernell2errorintro}.

\subsection{Storage requirements of $\calH^2$-farfield approximations}

As the final result of this section we discuss the following estimate on the storage requirements of the farfield of variable-order $\calH^2$-approximations.
\begin{lemma}\label{lem:H2complexity}
Let $g\in V^{\calH}$ with $\alpha\in\bbN_0$ and  $\beta\in\bbN$. Then the storage requirements for the coefficients of all leafs $t\in L_{I\times I}^+$ are bounded by
\[
C_{\calH^2}((\alpha+\beta)^{\delta d}|I|),
\]
i.e., they are linear with respect to the cardinality of the underlying index set $I$. The constant $C_{\calH^2}$ is independent of the depth of $T_{I\times I}$ and depends only on $\delta$, $d$, $C_{\text{sp}}$, and the shape of $T_I$ (see \cref{sec:H2appendix} for a precise statement).
\end{lemma}
\begin{proof}
We use the framework provided in \cite[Chapter 3.8]{Bor2010}. \Cref{lem:boundedrankdistr} shows that the rank as given by \cref{eq:rankdistribution} yields a $(1,\alpha,\beta,\delta d,C_{\text{ab}})$-bounded rank distribution in the sense of \cite[Definition 3.44]{Bor2010}, see also \cref{def:boundedrankdistr}. \Cref{lem:regularclustertree} yields that $T_I$ is a $(C_{rc},\alpha,\beta,\delta d,C_{\text{ab}})$-regular cluster tree in the sense of \cite[Definition 3.47]{Bor2010}, see also \cref{def:regularclustertree}, with $C_{rc}$ given as in \cref{eq:crc}. The assertion follows from \cite[Corollary 3.49]{Bor2010}, see also \cref{lem:complexitylemma}.
\end{proof}

\section{Proofs for the $\calH^2$-sample covariance estimator}\label{sec:MC}

\subsection{Algorithmic realization of $\Pi^\calH$ applied to simple tensors}\label{sec:Pih}

We first discuss that \cref{alg:Pih} indeed implements $\Pi^\calH(z_h\otimes z_h)$, $z_h\in V_h$, and recall that $\Pi^\calH$ is the $L^2$-projection onto $V^\calH$, cf.\ \cref{eq:defPih}. The discussion is split into two parts. First, we justify lines 4,7, and 8, which are responsible for the mathematical correctness of the algorithm. Lines 5 and 6 are responsible for computational efficiency.

\subsubsection{Lines 4,7, and 8}
We start by observing that \cref{eq:L2local} implies that for any $z_h\in V_h$ we have
\[
\Pi^\calH(z_h\otimes z_h)=\sum_{t\times s\in L_{I\times I}^+}\Pi_{t\times s}^\calH (z_h|_{D_t}\otimes z_h|_{D_s})+\sum_{t\times s\in L_{I\times I}^-}z_h|_{D_t}\otimes z_h|_{D_s},
\]
where $\Pi_{t\times s}^\calH (z_h|_{D_t}\otimes z_h|_{D_s})=u_{t\times s}^{\text{pw}}\in\calP_{t\times s}^{\text{pw}}$ are the solutions of the local variational problems
\begin{align}\label{eq:localL2}
\begin{aligned}
&\text{Find $u_{t\times s}^{\text{pw}}\in\calP_{t\times s}^{\text{pw}}$ s.t.}\\
&\qquad\qquad\qquad\text{$(u_{t\times s}^{\text{pw}},q_{t\times s}^{\text{pw}})_{L^2(D_t\times D_s)}=(z_h|_{D_t}\otimes z_h|_{D_s},q_{t\times s}^{\text{pw}})_{L^2(D_t\times D_s)}$}\\
&\text{for all $q_{t\times s}^{\text{pw}}\in\calP_{t\times s}^{\text{pw}}$,}
\end{aligned}
\end{align}
for all $t\times s\in L_{t\times s}^+$. 

To further simplify the variational problem we note that $\calP_{t\times s}^{\text{pw}}$, cf.\ \cref{eq:Ptspw}, inherits the tensor product structure of $\calP_{t\times s}$, cf.\ \cref{eq:Pts}, meaning that
\[
\calP_{t\times s}^{\text{pw}}
=
\calP_t^{\text{pw}}\otimes\calP_s^{\text{pw}}
\]
for all $t\times s\in L_{I\times I}^+$, where
\[
\calP_t^{\text{pw}} = \{f\in L^2(D_t)\colon f=\calI_{t_0}^tq,t_0\in L_t,q\in\calP_{k_{p-\level(t)}}\big|_{D_t}\},
\]
for all $t\in T_I$. Thus, \cref{eq:localL2} is equivalent to solving the finite dimensional variational problems
\[
\text{
Find $u_r^{\text{pw}}\in\calP_{t}^{\text{pw}}$ s.t.~$(u_r^{\text{pw}},q_r^{\text{pw}})_{L^2(D_r)}=(z_h|_{D_r},q_r^{\text{pw}})_{L^2(D_r)}$ for all $q_r^{\text{pw}}\in\calP_r^{\text{pw}}$,}
\]
for $r\in\{t,s\}$ and setting $u_{t\times s}^{\text{pw}}=u_t^{\text{pw}}\otimes u_s^{\text{pw}}$. Fixing appropriate nodal bases $\calP_r^{\text{pw}}=\operatorname{span}\{\psi_i^r\}_{i=1}^{m}$ with $m$ as in \cref{eq:rankdistribution} this is equivalent to solving the systems of linear equations
\begin{align}\label{eq:tensorlocalsle}
\bfQ_r\bfu_r=\bfq_r^h
\end{align}
with
\begin{align}\label{eq:tensorlocalL2matrices}
\bfQ_r=\big[(\psi_i^r,\psi_j^r)_{L^2(D_r)}\big]_{i,j=1}^m,
\quad
\bfq_r^h=\big[(z_h|_{D_r},\psi_i^r)_{L^2(D_r)}\big]_{i=1}^m,
\quad
\bfu_r=\big[\psi_i^r\big]_{i=1}^m,
\end{align}
for $r\in\{t,s\}$. The expression for $\bfq_r^h$ can be further simplified to
\begin{align*}
\bfq_r^h=\bfM_r\bfz_r^h,
\end{align*}
where $\bfM_r=\big[(\psi_i^r,\phi_j^r)_{L^2(D_r)}\big]_{i,j}$, $\psi_i^r\in\calP_r^{\text{pw}}$, $\phi_j^r\in V_j|_{D_r}$, is the moment matrix on $r$ and $\bfz_r^h$ is the coefficient vector of $z_h|_{D_r}$. We note that $\calP_t^{\text{pw}}=\calP_t$ for all $t\in L_I$, implying that
\begin{align*}
\bfM_r&=\big[(\psi_i^r,\phi_j^r)_{L^2(D_r)}\big]_{i,j},&&\psi_i^r\in\calP_r, \phi_j^r\in V_j|_{D_r},\\
\bfQ_r&=\big[(\psi_i^r,\psi_j^r)_{L^2(D_r)}\big]_{i,j=1}^m,&&\psi_i^r\in\calP_r.
\end{align*}
This justifies lines 4, 7, and 8 in \cref{alg:Pih}.

\subsubsection{Lines 5 and 6}\label{sec:H2SCEl56}
For practical purposes, it is preferable to avoid implementation of the spaces $\calP_r^{\text{pw}}$ due to computational cost and implementation complexity. Fortunately, it is known that vectors $\{\bfq_r^h\}_{r\in T_I}$ can be efficiently computed from $\{\bfq_r^h\}_{r\in L_I}=\{\bfM_r\bfz_r^h\}_{r\in L_I}$ and the transfer matrices $\{\bfE_t\}_{t\in T_I\setminus\{I\}}$ using the $\calH^2$-forward transformation algorithm \cite{Bor2010}. This is done in line 5.

For the computation of the matrices $\{\bfQ_r\}_{r\in T_I}$ one can apply the same iterative scheme, but in both, ansatz and test function. This is done in line 6.

\subsection{Computational complexity of \cref{alg:Pih}}\label{sec:Pihcompl}

We show the following theorem.
\begin{theorem}\label{thm:samplelinearcomplexity}
There is a constant $C>0$, independent of $\alpha,\beta,\delta,d$, and $|I|$, such that \cref{alg:Pih} requires at most $C(\alpha+\beta)^{2\delta d}|I|$ operations.
\end{theorem}
\begin{proof} We go line by line through the algorithm.
\begin{description}
\item[Lines 4 and 5:] Line 4 and 5 correspond to the forward transformation. It is a classical result from the literature, see \cite[Lemma 3.45 and 3.48]{Bor2010}, using the same constants as in the proof of \cref{lem:H2complexity}, that this operation can be done in at most $C_{\calH^2}(\alpha+\beta)^{\delta d}|I|$ operations.
\item[Line 6:] Using the Kronecker representation of the matrix products in line 6, cf.\ \cref{rem:algorem}, the same argumentation implies that line 6 can be done in at most $2C_{\calH^2}(\alpha+\beta)^{2\delta d}|I|$ operations.
\item[line 7:] Solving the local systems $t\in T_I$ is achievable in at most $3K_t^3$ complexity if a dense solver is used, with $K_t$ given as in \cref{eq:rankdistribution}. \cite[Lemma 3.45 and 3.48]{Bor2010} with the same constants as in the proof of \cref{lem:H2complexity} yields that solving all local systems requires $3C_{\calH^2}(\alpha+\beta)^{2\delta d}|I|$ operations in total.
\item[Line 8:] Computing $\bfu_t\otimes\bfu_s$, $t\times s\in L_{I\times I}^+$, requires $K_tK_s$ operations. \cite[Lemma 3.49]{Bor2010} yields that this can be achieved in $C_{\calH_ 2}(\alpha+\beta)^{2\delta d}|I|$ operations. 
\end{description} 
This yields the assertion.
\end{proof}

\subsection{Proof of \cref{lem:covHL2errorintro}}
We start by analyzing the approximation error of covariance functions in $\Pi^\calH(V_h\otimes V_h)$, which will be helpful later on.
\begin{corollary}\label{cor:PihPimixerr}
	Let $2\theta<0$ and let $g\in H^\gamma(D\times D)$, $\gamma>0$, be $G^\delta(C_G,A,2\theta)$-asymptotically smooth. Choose $\alpha\in\bbN$ such that $\zeta^{-2\theta}\tilde{\rho}^\alpha<1$ with $\tilde{\rho}$ as in \cref{eq:rhotilde}. Then there is $\beta_0\in\bbN$ such that
	\[
	\big\|g-\Pi^\calH\Pi_h^{\mix} g\big\|_{L^2(D\times D)}
	\leq
	\frac{C_{\text{lc}}C_{\text{sp}}h_{\calH}^{-2\theta}\tilde{\rho}^{\beta}}{1-\zeta^{-2\theta}\tilde{\rho}^{\alpha}}+C_{L^2}^\otimes h^\gamma\|g\|_{H^\gamma(D\times D)}
	\]
	for all $\beta\geq\beta_0$ with $\tilde{\rho}$ as in \cref{eq:rhotilde}.
\end{corollary}
\begin{proof}
	Follows from stability of the $L^2$-projection,
	\[
	\big\|g-\Pi^\calH\Pi_h^{\mix} g\big\|_{L^2(D\times D)}
	\leq
	\big\|g-\Pi^\calH g\big\|_{L^2(D\times D)}
	+
	\big\|g-\Pi_h^{\mix} g\big\|_{L^2(D\times D)},
	\]
	\cref{eq:approximationestimate}, and \cref{thm:L2interpolationestimate}.
\end{proof}

This allows us to estimate the $\calH^2$-SCE by means of a standard argument.

\begin{theorem}\label{lem:covHL2error}
Let $2\theta<0$ and let $g\in H^\gamma(D\times D)$ be $G^\delta(C_G,A,2\theta)$-asymptotically smooth. Choose $\alpha\in\bbN$ such that $\zeta^{-2\theta}\tilde{\rho}^\alpha<1$ with $\tilde{\rho}$ as in \cref{eq:rhotilde}. Then there is $\beta_0\in\bbN$ such that
\[
\big\|g-E^{MC}[\Pi^\calH\Pi_h^{\mix}g]\big\|_{L_{\bbP}^2(\Omega,L^2(D\times D))}
\leq
\frac{C_{\text{lc}}C_{\text{sp}}h_{\calH}^{-2\theta}\tilde{\rho}^{\beta}}{1-\zeta^{-2\theta}\tilde{\rho}^{\alpha}}+\bigg(C_{L^2}^\otimes h^{\gamma}+ \frac{1}{\sqrt{M}}\bigg)\|g\|_{H^\gamma(D\times D)}
\]
for all $\beta\geq\beta_0$ with $\tilde{\rho}$ as in \cref{eq:rhotilde}.
\end{theorem}
\begin{proof}
We repeat the argument of \cite[Theorem 4.3]{BSZ2011}. To this end, we estimate
\begin{align*}
&\big\|g-E^{MC}[\Pi^\calH\Pi_h^{\mix}g]\big\|_{L_{\bbP}^2(\Omega,L^2(D\times D))}\\
&\qquad\leq
\big\|g-\Pi^\calH\Pi_h^{\mix}g\big\|_{L_{\bbP}^2(\Omega,L^2(D\times D))}
+
\big\|\Pi^\calH\Pi_h^{\mix}g-E^{MC}[\Pi^\calH\Pi_h^{\mix}g]\big\|_{L_{\bbP}^2(\Omega,L^2(D\times D))}\\
&\qquad\leq
\big\|g-\Pi^\calH\Pi_h^{\mix}g\big\|_{L^2(D\times D)}
+
\big\|\Pi^\calH\Pi_h^{\mix}g-E^{MC}[\Pi^\calH\Pi_h^{\mix}g]\big\|_{L_{\bbP}^2(\Omega,L^2(D\times D))}.
\end{align*}
The first term is estimated by \cref{cor:PihPimixerr}. For the second term we note that $\Pi^\calH(V_h\otimes V_h)$ is a closed subspace of $L^2(D\times D)$, and thus a Hilbert space when equipped with the $L^2(D\times D)$ inner product. Thus, the second term is estimated by the standard Monte Carlo estimate
\[
\big\|\Pi^\calH\Pi_h^{\mix}g-E^{MC}[\Pi^\calH\Pi_h^{\mix}g]\big\|_{L_{\bbP}^2(\Omega,L^2(D\times D))}
\leq
\frac{1}{\sqrt{M}}\big\|\Pi^\calH\Pi_h^{\mix}g\big\|_{L_{\bbP}^2(\Omega,L^2(D\times D))},
\]
see \cite[Lemma 4.1]{BSZ2011} for a proof. $L^2$-stability of $\Pi^h\Pi_h^{\mix}$ implies
\[
\frac{1}{\sqrt{M}}\big\|\Pi^\calH\Pi_h^{\mix}g\big\|_{L_{\bbP}^2(\Omega,L^2(D\times D))}
\leq
\frac{1}{\sqrt{M}}\|g\|_{L_{\bbP}^2(\Omega,L^2(D\times D))},
\]
which yields the assertion.
\end{proof}
Thus, we have shown \cref{lem:covHL2errorintro}.

\section{Derivation of \cref{alg:MLreduction} and complexity analysis}\label{sec:MLMCalgo}
\subsection{Principal idea of the algorithm}
Reformulating the sum over the levels in the $\calH^2$-MLSCE, we seek an efficient implementation of the \emph{multilevel reduction}
\begin{align}\label{eq:reduction}
\Pi^{\calH_L}\colon\bigg[\bigtimes_{\ell=0}^LW_\ell^{\calH}\bigg]\to W_L^{\calH},\qquad\big[g^{\calH_\ell}\big]_{\ell=0}^L\mapsto \tilde{g}^{\calH_L}=\Pi^{\calH_L}\sum_{\ell=0}^Lg^{\calH_\ell},
\end{align}
with $W^{\calH}_\ell$ as in \cref{eq:W}. Exploiting \cref{eq:L2local,eq:adminimpladmin}, the sum in \cref{eq:reduction} can be rewritten as
\begin{align*}
&\Pi^{\calH_L}\sum_{\ell=0}^Lg^{\calH_\ell}\\
&\qquad=
\sum_{\ell=0}^L\Pi^{\calH_L}\bigg(\sum_{t\times s\in L_{I_\ell\times I_\ell}}g^{\calH_\ell}|_{D_t\times D_s}\bigg)\\
&\qquad=
\sum_{\ell=0}^L\Pi^{\calH_L}\bigg(\sum_{t\times s\in L_{I_\ell\times I_\ell}^+}g^{\calH_\ell}|_{D_t\times D_s}+\sum_{t\times s\in L_{I_\ell\times I_\ell}^-}g^{\calH_\ell}|_{D_t\times D_s}\bigg)\\
&\qquad=
\sum_{\ell=0}^L\left(
\sum_{t\times s\in L_{I_\ell\times I_\ell}^+}\Pi_{t\times s}^{\calH_L}g^{\calH_\ell}|_{D_t\times D_s}\phantom{\sum_{\substack{t'\times s'\in L_{I_L\times I_L}^-\\t'\times s'\subset t\times s}}}\right.\\
&\phantom{\qquad=\sum_{\ell=0}^L\bigg(}\left.+
\sum_{t\times s\in L_{I_\ell\times I_\ell}^-}\left(
\sum_{\substack{t'\times s'\in L_{I_L\times I_L}^+\\t'\times s'\subset t\times s}}\Pi_{t'\times s'}^{\calH_L}g^{\calH_\ell}|_{D_t\times D_s}
+
\sum_{\substack{t'\times s'\in L_{I_L\times I_L}^-\\t'\times s'\subset t\times s}}\Big(g^{\calH_\ell}|_{D_t\times D_s}\Big)\Big|_{t'\times s'}
\right)
\right).
\end{align*}
\Cref{alg:MLreduction} implements this sum in a recursive way. To this end, the details of computing the entities $\Pi_{t\times s}^{\calH_L}g^{\calH_\ell}|_{D_t\times D_s}$ in the sum over $t\times s\in L_{I_\ell\times I_\ell}^+$ are discussed in \cref{sec:admtoadm}. For the second sum, the special case when $t\times s\in L_{I_\ell\times I_\ell}^-$ with $t\times s\in L_{I_{\ell+1}\times I_{\ell+1}}^+$ is discussed in \cref{sec:inadmtoadm}. This explains line 16 of \cref{alg:MLreduction}. The remaining cases of the second sum, $t\times s\in L_{I_\ell\times I_\ell}^-$ with $t\times s\not\in L_{I_{\ell+1}\times I_{\ell+1}}^+$ are dealt with recursively, in lines 21--26 in \cref{alg:MLreduction}. The remaining parts of \cref{alg:MLreduction} are required for efficiency, and are discussed during the complexity analysis in \cref{sec:mlredcomplexity}.

\subsection{Projecting admissible block-clusters to admissible block-clusters}\label{sec:admtoadm}
We consider the case where $t\times s\in L_{I_\ell\times I_\ell}^+$ and $t\times s\in L_{I_L\times I_L}^+$, i.e., $t\times s$ is an admissible block-cluster in both block-cluster trees. For these block-clusters, computing $\Pi^{\calH_L}_{t\times s}g^{\calH_\ell}|_{D_t\times D_s}\in\calP_{t\times s}^{\text{pw},L}$, $g^{\calH_\ell}\in W^{\calH}_\ell$, amounts to the solution of
\begin{align*}
&\text{	Find $\Pi^{\calH_L}_{t\times s}g^{\calH_L}|_{D_t\times D_s}\in\calP_{t\times s}^{\text{pw},L}$ s.t.}\\
&\qquad\qquad\qquad\text{$\Big(\Pi^{\calH_L}_{t\times s}g^{\calH_L}|_{D_t\times D_s},q_{t\times s}^{\text{pw},L}\Big)_{L^2(D_t\times D_s)}=(g^{\calH_\ell}|_{D_t\times D_s},q_{t\times s}^{\text{pw},L})_{L^2(D_t\times D_s)}$}\qquad\qquad\qquad\\
&\text{for all $q_{t\times s}^{\text{pw},L}\in\calP_{t\times s}^{\text{pw},L}$.}
\end{align*}
This is a finite dimensional variational problem which can be written as
\begin{align}\label{eq:mlmcfarfieldreduction}
\bfQ_{t}\bfu_{t\times s}^{(L)}\bfQ_{s}^\intercal=\bfR_{t}^{(L,\ell)}\bfU_{t\times s}^{(\ell)}\Big(\bfR_{s}^{(L,\ell)}\Big)^\intercal,
\end{align}
with $\bfQ_r$, $r\in\{s,t\}$, as in \cref{eq:tensorlocalL2matrices}, $\bfu_{t\times s}^{(L)}$ and $\bfu_{t\times s}^{(\ell)}$ the coefficient matrices of $g^{\calH_L}|_{D_t\times D_s}$ and $g^{\calH_\ell}|_{D_t\times D_s}$, and 
\[
\bfR_r^{(L,\ell)}=\Big[\Big(\psi_i^{(r,L)},\psi_j^{(r,\ell)}\Big)_{L^2(D_r)}\Big]_{\substack{i=1,\ldots,K_r^{(L)},\\j=1,\ldots,K_r^{(\ell)}}}\in\bbR^{K_r^{(L)}\times K_r^{(\ell)}},
\]
for all $\psi_i^{(r,L)}\in\calP_r^{\text{pw},L}$ and $\psi_i^{(r,\ell)}\in\calP_r^{\text{pw},\ell}$, $r\in\{s,t\}$.

\subsection{Projecting inadmissible leaf block-clusters to admissible block-clusters}\label{sec:inadmtoadm}
We consider the case $t\times s\in L_{I_\ell\times I_\ell}^-$ and $t\times s\in L_{I_L\times I_L}^+$. Upon noting that it holds $g^{\calH_\ell}|_{D_t\times D_s}\in V_{h_\ell}|_{D_s}\otimes V_{h_\ell}|_{D_t}$ for all $g^{\calH_\ell}\in W^{\calH}_\ell$ we readily remark that
\begin{align*}
&\text{	Find $\Pi^{\calH_L}_{t\times s}g^{\calH_L}|_{D_t\times D_s}\in\calP_{t\times s}^{\text{pw},L}$ s.t.}\\
&\qquad\qquad\qquad\text{$(\Pi^{\calH_L}_{t\times s}g^{\calH_L}|_{D_t\times D_s},q_{t\times s}^{\text{pw},L})_{L^2(D_t\times D_s)}=(g^{\calH_\ell}|_{D_t\times D_s},q_{t\times s}^{\text{pw},L})_{L^2(D_t\times D_s)}$}\qquad\qquad\qquad\\
&\text{for all $q_{t\times s}^{\text{pw},L}\in\calP_{t\times s}^{\text{pw},L}$,}
\end{align*}
is a finite dimensional variational problem which can be rewritten as
\begin{align}\label{eq:mlmcfarfieldreduction2}
\bfQ_{t}\bfu_{t\times s}^{(L)}\bfQ_{s}^\intercal=\bfN_{t}^{(L,\ell)}\bfu_{t\times s}^{(\ell)}\Big(\bfN_{s}^{(L,\ell)}\Big)^\intercal.
\end{align}
As in the previous subsection, $\bfu_{t\times s}^{(L)}$ and $\bfu_{t\times s}^{(\ell)}$ are the coefficient matrices of $g^{\calH_L}|_{D_t\times D_s}$ and $g^{\calH_\ell}|_{D_t\times D_s}$, and 
\[
\bfN_r^{(L,\ell)}=\Big[\Big(\psi_i^{(r,L)},\phi_j^{(r,\ell)}\Big)_{L^2(D_r)}\Big]_{\substack{i=1,\ldots,K_r^{(L)},\\j=1,\ldots,\dim(V_{h_\ell}|r)}}\in\bbR^{K_r^{(L)}\times \dim(V_{h_\ell}|_{D_r})},
\]
for all $\psi_i^{(r,L)}\in\calP_r^{\text{pw},L}$ and $\phi_i^{(r,\ell)}\in V_{h_\ell}|_{D_r}$, $r\in\{s,t\}$.

\subsection{Complexity of the multilevel $\calH^2$-reduction algorithm}\label{sec:mlredcomplexity}
The following theorem shows that \cref{alg:MLreduction} has linear complexity in $|I_L|$.
\begin{theorem}\label{lem:reductioncomplexity}
There is a constant $C_{\text{ML}}=C_{\text{ML}}(C_{\calH^2},C_{\min},C_{\text{ab}},C_{\text{uni}},n_{\min},\delta,d)$ such that the computational cost of \cref{alg:MLreduction} are bounded by
\[
C_{\text{ML}}\frac{(\alpha+\beta)^{\lceil 3\delta d\rceil}}{(\alpha+\beta)^{\delta d}}|I_L|.
\]
\end{theorem}
\begin{proof}
We first list the computational cost of every step.
\begin{description}
\item[Lines 4 and 5:] As discussed in \cref{sec:H2SCEl56}, both lines correspond the forward transformation, which has linear complexity.
\item[Line 6:] The computational efforts for each $t'\in\children(t)$ are bounded by
\[
C_{\min}n_{\min}\Big(K_t^{(L)}K_{t'}^{(L)}+C_{\min}n_{\min}K_{t'}^{(L)}\Big)\leq 2C_{\min}^2n_{\min}^2K_t^{(L)}K_{t'}^{(L)}.
\]
Asymptotically, this is cheaper than Lines 4 and 5, implying that line 6 has linear complexity as well.
\item[Line 8--10:] 	We first note that $T_{I_\ell}$ is a $(C_{rc},\alpha,\beta+(L-\ell)\alpha,\delta d,C_{\text{ab}})$-bounded as well as a $(C_{rc},\alpha,\beta,\delta d,C_{\text{ab}})$-regular cluster tree with $C_{\text{rc}}$ as in \cref{eq:crc}. \Cref{lem:complexitylemma} implies that $\{\bfR_{t}^{(L,\ell)}\}_{t\in T_{I_\ell}}$ in line 9 can be computed in
\[
\calO\bigg(\frac{(\alpha(L-\ell+1)+\beta)^{3\delta d}}{(\alpha+\beta)^{\delta d}}|I_\ell|\bigg)
\]
operations, with the constant being independent of $\ell$ and $L$. Following the same arguments as for line 6, computing $\{\bfN_{t}^{(L,\ell)}\}_{t\in T_{I_\ell}}$ in line 10 requires at most
\[
\calO\bigg(C_{\min}^2n_{\min}^2\frac{(\alpha(L-\ell+1)+\beta)^{3\delta d}}{(\alpha+\beta)^{\delta d}}|I_\ell|\bigg)
\]
operations, with the constant being independent of $\ell$ and $L$.
\item[Line 16:] The computational cost for solving \cref{eq:mlmcfarfieldreduction} are given by
\[
\sum_{r\in\{s,t\}}\bigg(2\Big(K_r^{(L)}\Big)^3+\Big(K_r^{(L)}\Big)\Big(K_r^{(\ell)}\Big)^2\bigg)
\leq
3\sum_{r\in\{s,t\}}\Big(K_r^{(L)}\Big)^3
\]
and arise for all $t\times s\in L_{I_\ell\times I_\ell}^+$, while the efforts for \cref{eq:mlmcfarfieldreduction2} are given by
\[
\sum_{r\in\{s,t\}}\bigg(2\Big(K_r^{(L)}\Big)^3+\Big(K_r^{(L)}\Big)C_{\min}^2n_{\min}^2\bigg)
\leq
3C_{\min}^2n_{\min}^2\sum_{r\in\{s,t\}}\Big(K_r^{(L)}\Big)^3
\]
and arise for all $t\times s\in L_{I_\ell\times I_\ell}^-\cap L_{I_L\times I_L}^+$.
\item[Line 19:] The prolongation from $V_{h_\ell}|_{D_t}\otimes V_{h_\ell}|_{D_s}$ to $V_{h_{\ell+1}}|_{D_t}\otimes V_{h_{\ell+1}}|_{D_s}$ can be accomplished in at most $2C_{\text{uni}}C_{\min}^3n_{\min}^3$ operations.
\item[Line 22:] For all $t'\times s'\in\children(t\times s)\cap L_{I_{\ell+1}\times I_{\ell+1}}^+$ we need to solve \cref{eq:mlmcfarfieldreduction2} on the level pair $(L,\ell+1)$ instead of $(L,\ell)$, i.e.,
\[
\bfQ_{t'}\bfu_{t'\times s'}^{(L)}\bfQ_{s'}^\intercal=\bfN_{t'}^{(L,\ell+1)}\bfu_{t'\times s'}^{(\ell+1)}\Big(\bfN_{s'}^{(L,\ell+1)}\Big)^\intercal.
\]
The cost for a given $t\times s\in L_{I_\ell\times I_\ell}^-\setminus L_{I_L\times I_L}^+$ are thus bounded by
\[
\sum_{t'\times s'\in\children(t\times s)}\sum_{r\in\{s',t'\}}\bigg(2\Big(K_r^{(L)}\Big)^3+\Big(K_r^{(L)}\Big)C_{\min}^2n_{\min}^2\bigg)
\leq
3C_{\min}^2C_{\text{ab}}^2n_{\min}^2\sum_{r\in\{s,t\}}\Big(K_r^{(L)}\Big)^3.
\]
\item[Overall cost:]
Considering that lines 8--11 and 14--29 are within a for-loop, taking the sum over all values of $\ell$ yields
\begin{align*}
	&\sum_{\ell=0}^{L-1}\calO\Big(C_{\min}^2C_{\text{ab}}^2n_{\min}^2\frac{(\alpha(L-\ell+1)+\beta)^{3\delta d}}{(\alpha+\beta)^{\delta d}}+C_{\text{uni}}C_{\min}^3n_{\min}^3\Big)|I_\ell|\\
	&\quad\leq |I_0|\sum_{\ell=0}^{L-1}\calO\Big(C_{\min}^2C_{\text{ab}}^2n_{\min}^2\frac{(\alpha(L-\ell+1)+\beta)^{3\delta d}}{(\alpha+\beta)^{\delta d}}+C_{\text{uni}}C_{\min}^3n_{\min}^3\Big)C_{\text{uni}}^{\ell}\\
	&\quad\leq |I_0|\calO\bigg(C_{\min}^2C_{\text{ab}}^2n_{\min}^2\sum_{\ell=0}^{L-1}\frac{(\alpha(L-\ell+1)+\beta)^{3\delta d}}{(\alpha+\beta)^{\delta d}}C_{\text{uni}}^{\ell}+C_{\text{uni}}C_{\min}^3n_{\min}^3\frac{C_{\text{uni}}^L-1}{C_{\text{uni}}-1}\bigg).
\end{align*}
We note that
\begin{align*}
	\sum_{\ell=0}^{L-1}\frac{(\alpha(L-\ell+1)+\beta)^{3\delta d}}{(\alpha+\beta)^{\delta d}}C_{\text{uni}}^{\ell}
	&\leq
	C_{\text{uni}}^{L}\sum_{\ell=0}^{L}\frac{((\alpha+\beta)+\alpha\ell)^{3\delta d}}{(\alpha+\beta)^{\delta d}}C_{\text{uni}}^{-\ell}\\
	&\leq
	\frac{C_{\text{uni}}^{L}}{(\alpha+\beta)^{\delta d}}\sum_{\ell=0}^{\infty}((\beta+\alpha)+\alpha\ell)^{\lceil 3\delta d\rceil}C_{\text{uni}}^{-\ell}
\end{align*}
where
\[
\sum_{\ell=0}^\infty(\beta+\alpha\ell)^kq^\ell\leq\bigg(1+\frac{1}{1-q}\bigg(\frac{q}{1-q}+\frac{1}{2}\bigg)^k k!\bigg)(\alpha+\beta)^k
\]
for all $q\in[0,1)$ and $k\in\bbN_0$ due to \cite[Lemma 3.50 and 3.51]{Bor2010}.
The assertion follows with $|I_0|C_{\text{uni}}^L=|I_L|$.
\end{description}
\end{proof}

\section{Numerical experiments}\label{sec:experiments}
\begin{figure}
\centering
\includegraphics[height=0.3\textwidth,clip=true,trim=340 120 330 0]{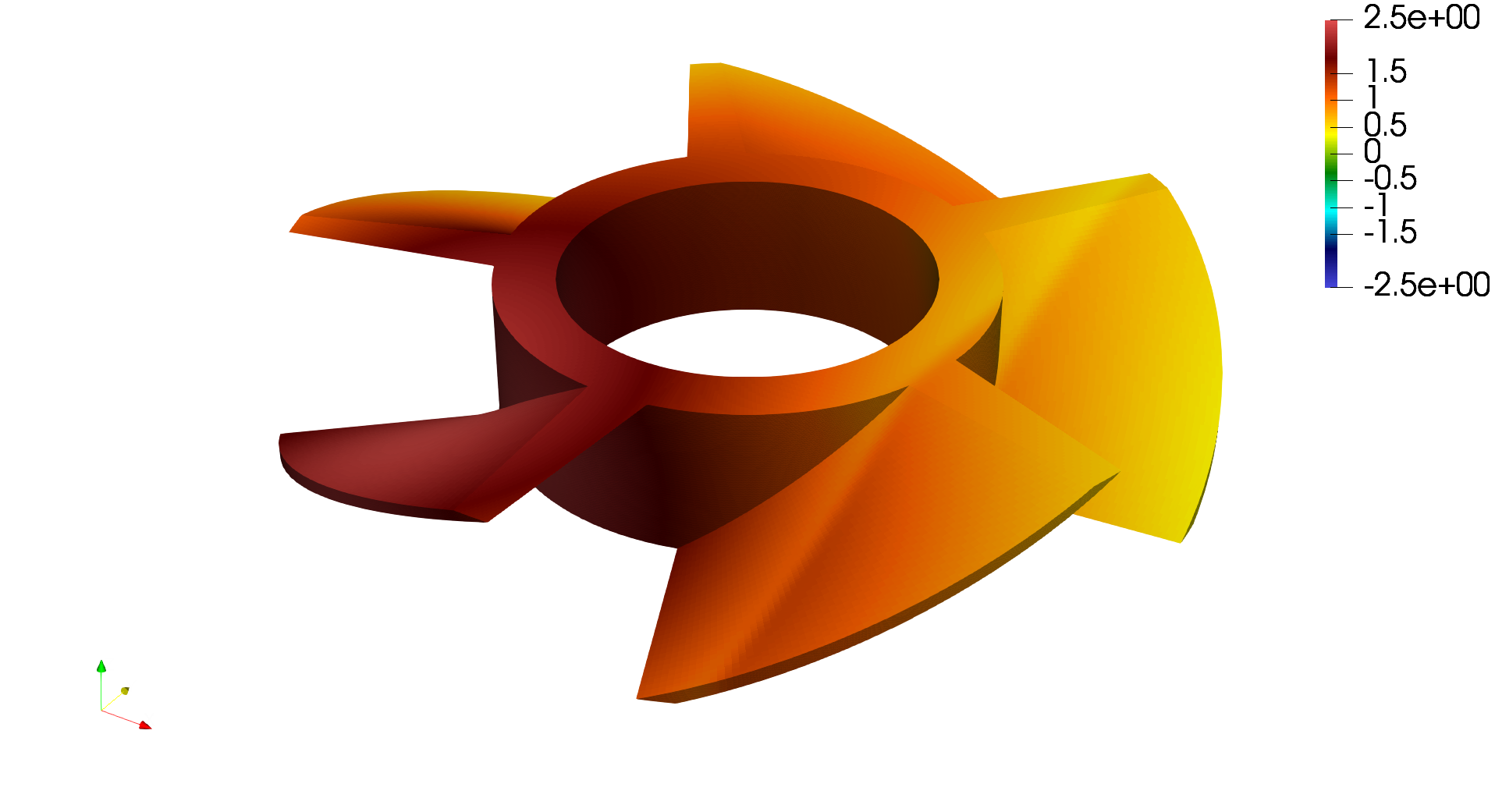}
\includegraphics[height=0.3\textwidth,clip=true,trim=340 120 330 0]{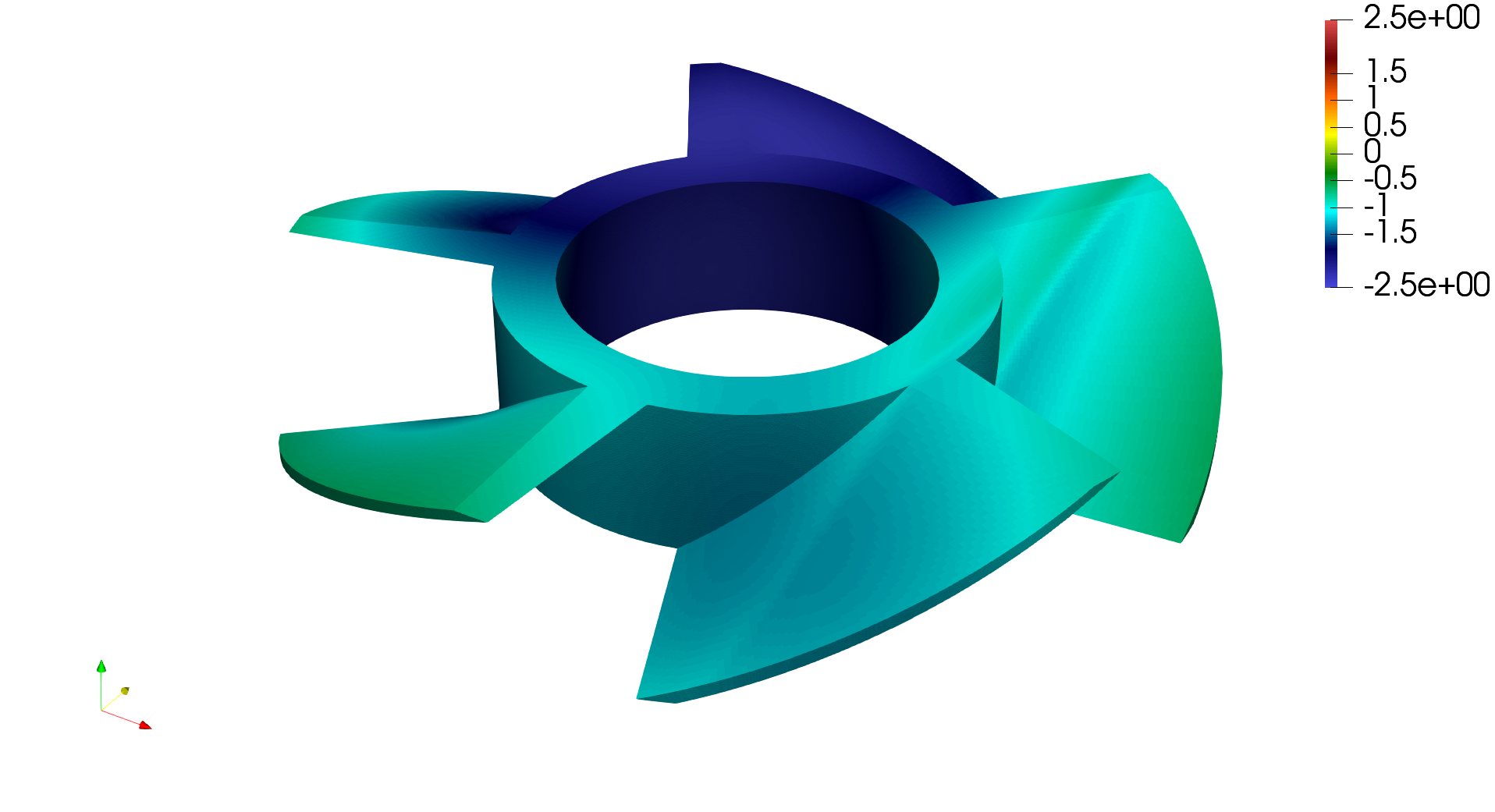}\\
\includegraphics[height=0.3\textwidth,clip=true,trim=340 120 330 0]{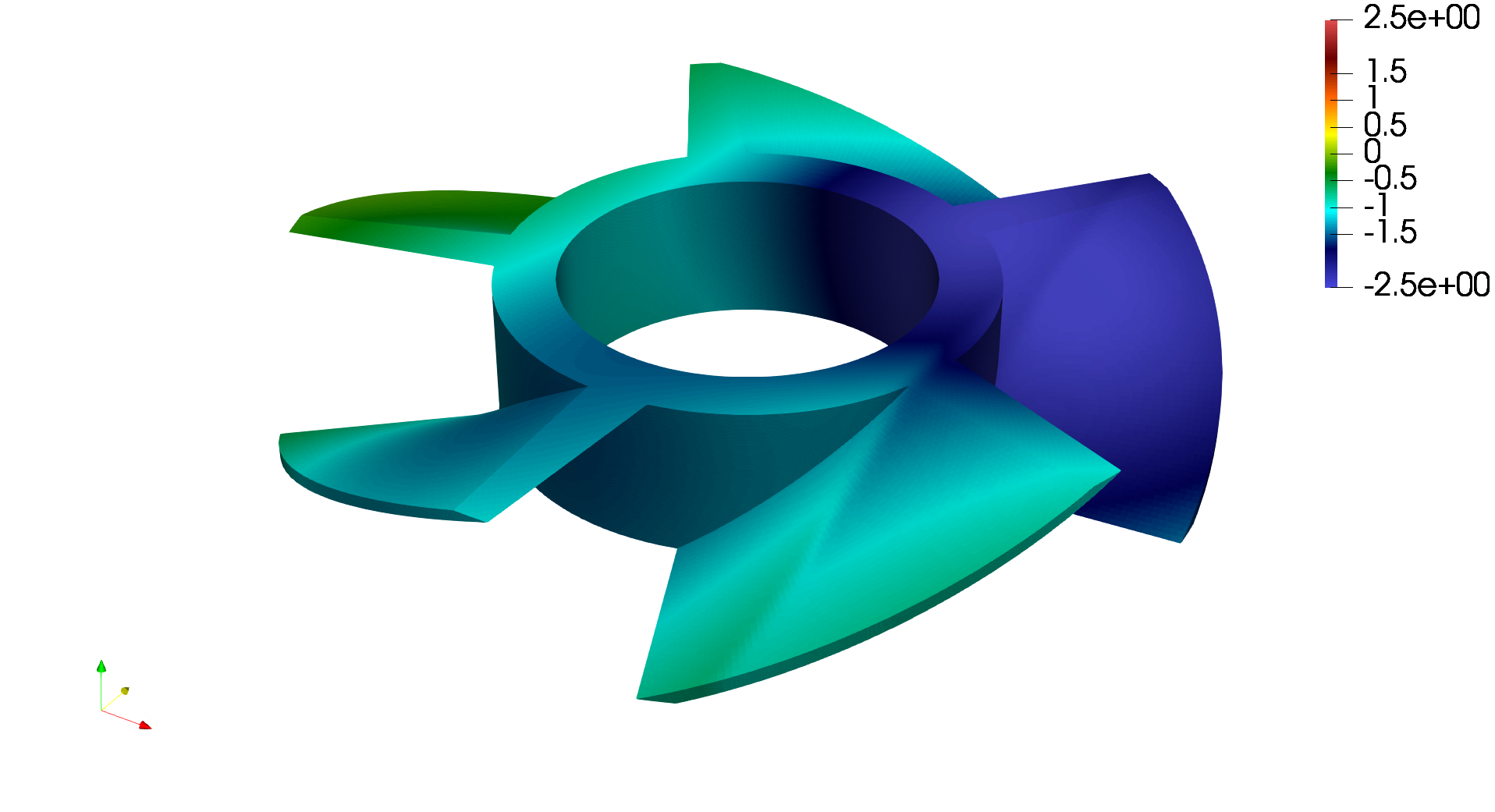}
\includegraphics[height=0.3\textwidth,clip=true,trim=340 120 0 0]{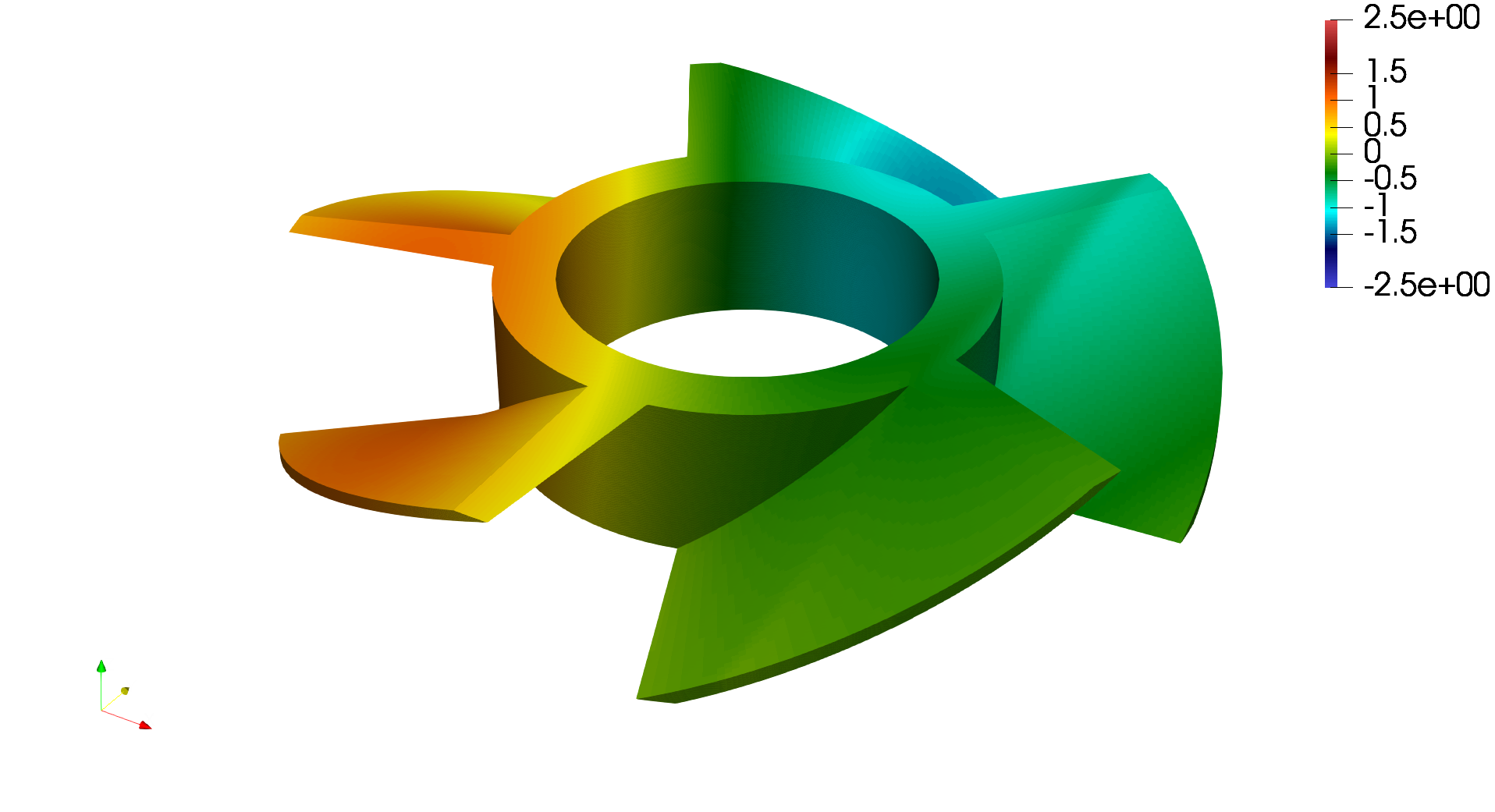}
\caption{\label{fig:turbine}Sample realizations of the centered Gaussian process with $G^{3/2}$-asymptotically smooth covariance function taken for the numerical experiments.}
\end{figure}

For our numerical experiments we aim at estimating the covariance of a Gaussian random field at the surface $\partial D$ of a turbine geometry $D$, see \cref{fig:turbine}, i.e., on a two-dimensional manifold embedded into $\bbR^3$. The radius of the turbine to the end of the blades is 1.5 and the whole geometry is contained in the box $[-1.5,1.5]\times[-1.5,1.5]\times[0,0.5]$.
To that end, we prescribe a reference Gaussian random field in terms of a Karhunen-Lo\'eve expansion, i.e.,
\[
\calZ(\omega,x)=\sum_{k=0}^\infty\sqrt{\lambda_k}\varphi(\bfx)Y_k(\omega),
\]
with $Y_k\sim U([-1,1])$ and $\{(\lambda_k,\varphi_k)\}_{k=0}^\infty$ the eigenpairs of the integral operator
\[
\calC\colon L^2(\partial D)\to L^2(\partial D),\qquad(\calC\varphi)(\bfx)=\int_{\partial D}g_\delta(\bfx,\bfy)\varphi(\bfy)\dd\sigma(\bfy).
\]
The covariance function $g_\delta$ is chosen as a modified Mat\'ern-$9/2$ kernel
\[
g_\delta(\bfx,\bfy)
=
\tilde{g}(\|\gamma_\delta(\bfx)-\gamma_\delta(\bfy)\|),
\qquad
\tilde{g}(r)
=
\bigg(
1+3r+\frac{27r^2}{7}+\frac{18r^3}{7}+\frac{27r^4}{35}
\bigg)\frac{e^{-3r}}{3},
\]
where
\[
\bfgamma_\delta\colon\partial D\to\bbR^3,\qquad
\bfgamma_\delta(x_1,x_2,x_3),
=
\begin{bmatrix}
(0.1 + (\delta-1)\Upsilon_\delta(2*x_1-1))x_1\\
x_2\\
x_3
\end{bmatrix}
\]
and
\[
\Upsilon_\delta(t) = \frac{\upsilon_\delta(1-t)}{\upsilon_\delta(1-t)+\upsilon_\delta(t)},
\qquad
\upsilon_\delta(t)=
\begin{cases}
0,&t\leq 0,\\
e^{-t^{\frac{1}{1-\delta}}},& t>0,
\end{cases}
\]
is a partition of unity of Gevrey class $\delta>1$ with $\Upsilon_\delta(t)=1$ for $t<0$ and $\Upsilon_\delta(t)=0$ for $t>1$, see, e.g., \cite{CR1996}. We set $\Upsilon_1\equiv 0$. Note that the $\delta-1$ factor in the definition of $\bfgamma_\delta$ vanishes for $\delta=1$, i.e., the analytic case, implying that the partition of unity disappears in this case. This makes the covariance function $g_\delta$ a $G^{\delta}$-asymptotically smooth kernel function for all $\delta\geq 1$. For our numerical experiments, we choose $\delta\in\{1,3/2,19/10\}$. Some illustrative samples for $\delta=3/2$ are shown in \cref{fig:turbine}, where the transition of the Gevrey partition of unity from zero to non-vanishing at $\{x_1=0.5\}$ is recognizable.

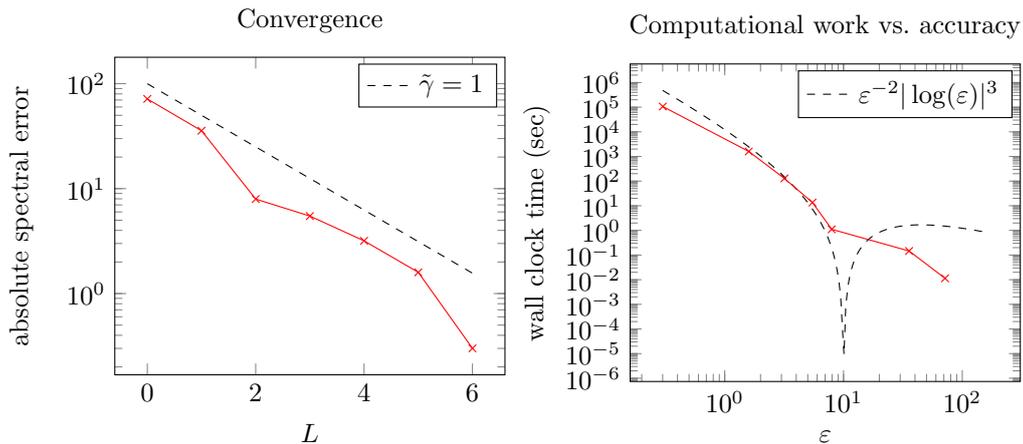
\begin{figure}
\centering
\begin{tikzpicture}
\begin{semilogyaxis}[
title=Convergence,
width=0.45\textwidth,
xlabel=$L$,
ylabel=absolute spectral error,
]
\addplot[domain=0:6,
samples=200,
dashed]{1e2*2^(-x)};
\addlegendentry{$\tilde{\gamma}=1$}
\addplot[red,
mark=x]table[x=L,y=eps]{data/MLMC_L6_refL6_CL1.000000_Gs1.000000.data};
\addlegendentry{$\delta=1$}
\addplot[blue,
mark=x]table[x=L,y=eps]{data/MLMC_L6_refL6_CL1.000000_Gs1.500000.data};
\addlegendentry{$\delta=3/2$}
\addplot[green,
mark=x]table[x=L,y=eps]{data/MLMC_L6_refL6_CL1.000000_Gs1.900000.data};
\addlegendentry{$\delta=19/10$}
\end{semilogyaxis}
\end{tikzpicture}
\begin{tikzpicture}
\begin{loglogaxis}[
title=Computational work vs.\ accuracy,
width=0.45\textwidth,
xlabel=$\varepsilon$,
ylabel=wall clock time (sec),
max space between ticks=10,
ymin=1e-4,
ymax=1e6,
xmax=1e2
]
\addplot[domain=1:71,
samples=200,
dashed]{15*(0.1*x)^(-2)*abs(ln(1/(0.1*x)))^3};
\addlegendentry{$\varepsilon^{-2}|\log(\varepsilon)|^3$}
\addplot[red,
mark=x]table[x=eps,y=time]{data/MLMC_L6_refL6_CL1.000000_Gs1.000000.data};
\addlegendentry{$\delta=1$}
\addplot[blue,
mark=x]table[x=eps,y=time]{data/MLMC_L6_refL6_CL1.000000_Gs1.500000.data};
\addlegendentry{$\delta=3/2$}
\addplot[green,
mark=x]table[x=eps,y=time]{data/MLMC_L6_refL6_CL1.000000_Gs1.900000.data};
\addlegendentry{$\delta=19/10$}
\end{loglogaxis}
\end{tikzpicture}
\caption{\label{fig:convergenceplot}Convergence plot of a realization of the $\calH^2$-MLSCE and corresponding computational work vs.\ accuracy with the sample numbers as in \cref{tab:samplenumbers}, cf.\ also \cref{thm:blfestimate} and \cref{thm:costcomplexity}.}
\end{figure}

The $\mathcal{H}^2$-implementation of the numerical experiments is based on the C++-Library Bembel \cite{DHK+2020}, with compression parameters $\alpha=1$, $\beta=2$, $\eta=0.8$, and $n_{\min}=4$. We choose piecewise constant finite element spaces $V_{h_\ell}=W_{h_\ell}$, $\ell=0,1,2,\ldots$, on uniformly refined quadrilateral meshes with $C_{\text{uni}}=4$ and $h_{\ell}\sim 2^{-\ell}$, leading to dimensions of the finite element spaces and covariance matrices as in \cref{tab:FEdimension}.
\begin{table}
\centering
\begin{tabular}{|c|ccccccc|}
	\hline
	$L$	& 0	& 1	& 2	& 3	& 4	& 5	& 6 \\\hline
	$\dim V_h = \dim W_h$ & 60 & 240 & 960 & 3 840 & 15 360 & 61 440 & 245 760 \\\hline
	$\dim (W_h\otimes W_h)$ & 3 600 & 57 600 & 921 600 & $\approx 14.7\cdot 10^6$ & $\approx 236\cdot 10^6$ & $\approx 3.77\cdot 10^9$ & $\approx 60.4\cdot 10^9$\\\hline
\end{tabular}
\caption{\label{tab:FEdimension} Dimensions of the used finite element spaces. The estimated covariance matrices are matrices in $\bbR^{\dim W_h\times\dim W_h}$, i.e., we have $\dim(W_h\otimes W_h)$ degrees of freedom.}
\end{table}
The Gaussian random field samples $\Pi_{h_\ell}\calZ$ are generated from a Karhunen-Lo\'eve expansion which is truncated at $10^{-3}h^\ell$ and computed from a pivoted Cholesky decomposition \cite{HPS2015}. According to \cref{thm:blfestimate} and \cref{thm:costcomplexity} it holds $\tilde{\gamma}=1$, and we can expect a linear convergence rate for our $\calH^2$-MLSCE, if the sample numbers are chosen proportional to \cref{thm:costcomplexity}. For our particular example we choose the sample numbers listed in \cref{tab:samplenumbers}.
\begin{table}
\centering
\begin{tabular}{|l|rrrrrrr|}
	\hline
	$L$	& 0	& 1	& 2	& 3	& 4	& 5	& 6 \\\hline
	$M_0$	& 1	& 4	& 64	& 576	& 4 096	& 25 600	& 147 456 \\
	$M_1$	& 	& 1	& 16	& 144	& 1 024	& 6 400	& 36 864 \\
	$M_2$	& 	& 	& 4	& 36	& 256	& 1 600	& 9 216 \\
	$M_3$	& 	& 	& 	& 9	& 64	& 400	& 2 304 \\
	$M_4$	& 	& 	& 	& 	& 16	& 100	& 576 \\
	$M_5$	& 	& 	& 	& 	& 	& 25	& 144 \\
	$M_6$	& 	& 	& 	& 	& 	& 	& 36 \\
	\hline
\end{tabular}
\caption{\label{tab:samplenumbers}Sample numbers chosen according to the case $2\tilde{\gamma}=d$ in \cref{thm:costcomplexity} for the numerical example.}
\end{table}
\Cref{fig:convergenceplot} shows that we reach indeed the predicted rate convergence rate of \cref{thm:MLestimatorerror} and a computational work vs.~accuracy as in \cref{thm:costcomplexity}. The spectral error was computed with a power iteration up to an absolute accuracy of $10^{-3}$. The computation times are measured in wall clock time and have been carried out in parallel with 64 threads on a compute server with 4TB RAM and two AMD EPYC 9354 32-Core Processors and hyperthreading disabled. The missing point for L=6 and $\delta=19/10$ is due to memory constraints.

\section{Conclusion}\label{sec:concl}
In this article, we considered the multilevel estimation of covariance functions which are $G^\delta$-asymp\-to\-ti\-cally smooth, $\delta\geq 1$. The naive approach to estimate the covariance function from discretized samples using the single level covariance estimator is computationally prohibitive due to the density of the arising covariance matrices and the slow convergence of the sample covariance estimator. To overcome these issues, we first generalized the classical $\calH^2$-approximation theory for asymptotically smooth kernels to Gevrey kernels. This allows to compress the arising covariance matrices by $\calH^2$-matrices in linear complexity with respect to the underlying approximation space. Secondly, we proposed and analyzed a \emph{$\calH^2$-formatted multilevel covariance sample estimator ($\calH^2$-MLCSE)}. This estimator exploits an approximate multilevel hierarchy in the $\calH^2$-approximation spaces to \emph{estimate the covariance in the same complexity as the mean}.

Alternatively to the approach proposed in this paper, a wavelet based method for estimating covariance functions was proposed in \cite{HHKS2021}. The advantage of such a wavelet method is that the wavelet-based approximation results also hold for finite smoothness of the covariance function, whereas the here presented $\calH^2$-approach requires asymptotically infinite smoothness. In contrast, the advantage of the $\calH^2$-approach in this paper is that no wavelet basis is required and that the presented algorithms can be integrated into the many readily available $\calH^2$-matrix codes.

\section*{Declarations}

\subsection*{Acknowledgement}
The author would like to express his sincere gratitude to Christoph Schwab for
the initial discussions on generalizing the $\calH^2$-matrix approximation theory to Gevrey kernels
and for critical and helpful comments during the writing of the manuscript.

\subsection*{Funding statement}
The author acknowledges the support from the DFG under Germany’s Excellence Strategy, project 390685813.

\subsection*{Conflicts of interest}
The authors declared that they have no conflict of interest.

\bibliographystyle{plain}
\bibliography{references}

\appendix

\section{Computation of $\calH^2$-related constants}\label{sec:H2appendix}
\begin{definition}[{\cite[Definition 3.44]{Bor2010}}]\label{def:boundedrankdistr}
	Let $T_I$ be a cluster tree and denote the number of interpolation points chosen in each cluster $t\in T_I$ by $K_t$. We say that $\{K_t\}_{t\in T_I}$ is a \emph{rank distribution}. We say that $\{K_t\}_{t\in T_I}$ is a \emph{$(C_{bn},\alpha,\beta,r,\xi)$-bounded rank distribution}, $C_{bn}\geq 1$, $\alpha>0$, $\beta\geq 0$, $r\geq 1$, $\xi\geq 1$, if
	\[
	\big|\big\{t\in T_I\colon K_t>(\alpha+\beta(\ell-1))^r\big\}\big|\leq C_{bn}\xi^{-\ell}|T_I|,\qquad\text{for all}~\ell\in\bbN.
	\]
\end{definition}
\begin{lemma}\label{lem:boundedrankdistr}
	Let $T_I$ be a cluster tree on the index set $I$ satisfying the assumptions of \cref{sec:ct}. Then $\{K_t\}_{t\in T_I}$ is a $(1,\alpha,\beta,\delta d,C_{\text{ab}})$-bounded rank distribution if the number of interpolation points in $(K_t)_{t\in T_I}$ are chosen according to \cref{eq:rankdistribution}, i.e.,
	\[
	K_t=\big\lceil(\beta+\alpha(p-\level(t)))^\delta\big\rceil^{d}
	\]
\end{lemma}
\begin{proof}
	The proof is analogy to the example in \cite[p.~64]{Bor2010}. We need to bound the number of clusters with
	\[
	K_t
	=
	\big\lceil(\beta+\alpha(p-\level(t)))^\delta\big\rceil^{d}
	\geq
	(\beta+\alpha(p-\level(t)))^{\delta d}
	>(\alpha+\beta(\ell -1))^{\delta d}.
	\]
	From this inequality we deduce that the clusters satisfying this constraint also satisfy $\level(t)<p+1-\ell$. Due to the assumptions in \cref{sec:ct} the number of such clusters is bounded by from above by $(C_{\text{ab}}^{p-\ell+2}-1)/(C_{\text{ab}}-1)$, and we obtain the assertion due to
	\[
	|T_I|
	\geq
	\frac{C_{\text{ab}}^{p+2}-1}{C_{\text{ab}}-1}
	=
	C_{\text{ab}}^{\ell}\frac{C_{\text{ab}}^{p-\ell+2}-C_{\text{ab}}^{-\ell}}{C_{\text{ab}}-1}
	\geq
	C_{\text{ab}}^{\ell}\frac{C_{\text{ab}}^{p-\ell+2}-1}{C_{\text{ab}}-1}.
	\]
\end{proof}
\begin{definition}[{\cite[Definitions 3.43 and 3.47]{Bor2010}}]\label{def:regularclustertree}
	Let $T_I$ be a cluster tree. We say that it is \emph{$(C_{rc},\alpha,\beta,r,\xi)$-bounded} with $C_{rc}\geq 1$, $\alpha>0$, $\beta\geq 0$, $r\geq 1$, $\xi>1$, if
	\begin{equation}\label{eq:rct4}
	\big|\big\{t\in L_I\colon|t|>(\beta+\alpha(\ell-1))^r\big\}\big|\leq C_{rc}\xi^{-\ell}|T_I|,\qquad\text{for all}~\ell\in\bbN,
	\end{equation}
	and
	\begin{align}
	|\children(t)|\leq C_{rc},\qquad\text{for all}~t\in T_I.\label{eq:rct1}
	\end{align}
	We say that $T_I$ is \emph{$(C_{rc},\alpha,\beta,r,\xi)$-regular}, if it is \emph{$(C_{rc},\alpha,\beta,r,\xi)$-bounded} and additionally satisfies
	\begin{align}
	|\children(t)|&\geq 2,&&\hspace*{-2cm}\text{for all}~t\in T_I\setminus L_I,\label{eq:rct2}\\
	(\alpha+\beta)^r&\leq C_{rc}|t|,&&\hspace*{-2cm}\text{for all}~t\in L_I.\label{eq:rct3}
	\end{align}
\end{definition}

\begin{lemma}\label{lem:regularclustertree}
	Let $T_I$ be a cluster tree with depth $p$ on the index set $I$ satisfying the assumptions of \cref{sec:ct}. Then $T_I$ is $(C_{rc},\alpha,\beta,\delta d,C_{\text{ab}})$-regular with
	\begin{equation}\label{eq:crc}
	C_{\text{rc}}=\max\bigg\{C_{\text{ab}},\frac{(\alpha+\beta)^{\delta d}}{n_{\min}},C_{\text{ab}}^{\frac{n_{\min}^{1/(\delta d)}-\beta+\alpha}{\alpha}+1}\bigg\}.
	\end{equation}
\end{lemma}
\begin{proof}
	\Cref{eq:ctnoc} implies  $2\leq|\children(t)|\leq C_{\text{ab}}$, $t\in T_I\setminus L_I$, which yields \eqref{eq:rct2} and \cref{eq:rct1} holds with $C_{\text{rc}}\geq C_{\text{ab}}$. Inserting the upper bound from \cref{eq:ctleafsize} into \cref{eq:rct3} yields
	\[
	\frac{(\alpha+\beta)^{\delta d}}{n_{\min}}\leq C_{\text{rc}}.
	\]
	
	Finally, the lower bound from \cref{eq:ctleafsize} implies that there are at most $C_{\text{ab}}^{p+1}$ leafs. The upper bound from \cref{eq:ctleafsize} and \cref{eq:rct4} with $\xi=C_{\text{ab}}$ then imply that $C_{rc}$ must satisfy
	\[
	C_{rc}
	\geq
	\begin{cases}
	\frac{C_{\text{ab}}^{p+\ell+1}}{|T_I|}&\text{for all}~\ell~\text{with}~(\beta+\alpha(\ell-1))^{\delta d}<n_{\min}\\
	0&\text{else}
	\end{cases}
	\]
	Solving $(\beta+\alpha(\ell-1))^{\delta d}<n_{\min}$ for $\ell$ implies $\ell<(n_{\min}^{1/(\delta d)}-\beta+\alpha)/\alpha$ which yields
	\[
	C_{\text{ab}}^{\frac{n_{\min}^{1/(\delta d)}-\beta+\alpha}{\alpha}+1}\leq C_{rc}
	\]
	due to $|T_I|\geq (C_{\text{ab}}^{p+2}-1)/(C_{\text{ab}}-1)$. Combining all conditions on $C_{rc}$ yields the assertion.
\end{proof}

\begin{lemma}[{\cite[Lemma 3.45]{Bor2010}}]\label{lem:borlem345}
	Let $T_I$ be a $(C_{\text{rc}},\alpha,\beta,r,\xi)$-bounded cluster tree and let $\{K_t\}_{t\in T_I}$ be a $(C_{\text{bn}},\alpha,\beta,r,\xi)$-bounded rank distribution. Define
	\begin{align}\label{eq:kt}
	k_t
	=
	\begin{cases}
	\max\{K_t,|t|\},&t\in L_I,\\
	\max\{K_t,\sum_{t'\in\children(t)}K_{t'}\},&t\in T_I\setminus L_I,
	\end{cases}
	\end{align}
	and $m\in\bbN$. Then there is a constant $C_{\text{cb}}=C_{\text{cb}}(C_{\text{rc}},C_{\text{bn}},r,\xi)\geq 1$ such that
	\[
	\sum_{t\in T_I}k_t^m\leq C_{\text{cb}}(\alpha+\beta)^{rm}|T_I|.
	\]
\end{lemma}
\begin{lemma}[{\cite[Lemma 3.48]{Bor2010}}]\label{lem:borlem348}
	Let $T_I$ be a $(C_{\text{rc}},\alpha,\beta,r,\xi)$-regular cluster tree. Then it holds
	\[
	|T_I|\leq \frac{2C_{\text{rc}}|I|}{(\alpha+\beta)^r}
	\]
\end{lemma}

\begin{lemma}[{Modification of \cite[Corollary 3.49]{Bor2010}}]\label{lem:complexitylemma}
	Let $T_I$ be $(C_{\text{rc}},\alpha,\beta,r,\xi)$-bounded and $(K_t)_{t\in T_I}$ be a $(C_{\text{bn}},\alpha,\beta,\xi)$-bounded rank distribution. Let $T_I$ be $(C_{\text{rc}},\alpha',\beta',r,\xi)$-regular and $T_{I\times I}$ be a block-cluster tree with sparsity constant $C_{\text{sp}}$. For $m\in\bbN$ and $\{k_t\}_{t\in T_I}$ defined as in \cref{eq:kt} it holds
	\[
	\sum_{t\in T_I}k_t^m\leq C_{\calH^2}\frac{(\alpha+\beta)^{rm}}{(\alpha'+\beta')^r}|I|
	\]
	with $C_{\calH^2}=2C_{\text{rc}}C_{\text{cb}}$.
\end{lemma}
\begin{proof}
	Combine \cref{lem:borlem345} and \cref{lem:borlem348}.
\end{proof}

\end{document}